\documentclass[11 pt,english]{article}
\usepackage[activeacute]{babel}
\usepackage{amsfonts} 
\usepackage{mathrsfs} 
\usepackage{amssymb, amsmath, amsfonts}
\usepackage{makeidx}
\usepackage{epsfig}
\usepackage{color}
\usepackage[T1]{fontenc}
\usepackage{graphics,graphicx,graphpap}
\usepackage[ansinew]{inputenc}
\usepackage{amsthm}

\newcommand\blfootnote[1]{%
  \begingroup
  \renewcommand\thefootnote{}\footnote{#1}%
  \addtocounter{footnote}{-1}%
  \endgroup
}

\usepackage{tikz}
\usepackage{subfigure}
\usetikzlibrary{calc, positioning, through, math, arrows, backgrounds}
\tikzstyle{vtx}=[inner sep=1pt,draw, shape=circle, font=\tiny]
\tikzstyle{line}=[inner sep=3pt,draw, shape=rectangle, line width = 3pt]
\tikzset{>=stealth}
\tikzstyle{lbl}=[inner sep = 1 pt, fill = white, midway]

\if@twoside \oddsidemargin 6pt \evensidemargin 6pt \marginparwidth 20pt
 \else \oddsidemargin 6pt \evensidemargin 6pt \marginparwidth 20pt \fi
\newcommand{\ZZ}{\mathbb{Z}}

\newcommand{\Aut}{\mathrm{Aut}}
\newcommand{\G}{\Gamma}
\renewcommand{\O}{\mathcal{O}}

\newcommand{\cC}{\mathcal{C}}

\newcommand{\cI}{\mathcal{I}}
\newcommand{\cA}{\mathcal{VT}}
\newcommand{\cXa}{\mathcal{X}_a}

\newcommand{\la}{\langle}
\newcommand{\ra}{\rangle}

\newcommand{\HTG}{\mathop{\rm HTG}}

\newcommand{\No}[1]{{#1}^*}
\newcommand{\vv}[1]{u_{#1}}

\newtheorem{theorem}{Theorem}[section]
\newtheorem{proposition}[theorem]{Proposition}
\newtheorem{corollary}[theorem]{Corollary}
\newtheorem{lemma}[theorem]{Lemma}

\theoremstyle{definition}

\newtheorem{problem}[theorem]{Problem}
\newtheorem{construction}[theorem]{Construction}
\newtheorem{assumption}[theorem]{Assumption}

\begin{document}

\begin{center}
\Large{\textbf{Cubic factor-invariant graphs of cycle quotient type - the alternating case}} \\ [+4ex]
\end{center}

\begin{center}
Brian Alspach{\small $^{a}$},\   
Primo\v z \v Sparl{\small $^{b,c,d,*}$}
\\

\medskip
{\it {\small
$^a$School of Information and Physical Sciences, University of Newcastle, Callaghan, NSW 2308, Australia\\
$^b$Institute of Mathematics, Physics and Mechanics, Ljubljana, Slovenia\\
$^c$University of Ljubljana, Faculty of Education, Ljubljana, Slovenia\\
$^d$University of Primorska, Institute Andrej Maru\v si\v c, Koper, Slovenia\\
}}
\end{center}

\blfootnote{
Email addresses: 
brian.alspach@newcastle.edu.au (B.~Alspach), primoz.sparl@pef.uni-lj.si (P. \v Sparl)
\\
* - corresponding author
}


\hrule

\begin{abstract}
We investigate connected cubic vertex-transitive graphs whose edge sets admit a partition into a $2$-factor $\cC$ and a $1$-factor that is invariant under a vertex-transitive subgroup of the automorphism group of the graph and where the quotient graph with respect to $\cC$ is a cycle. There are two essentially different types of such cubic graphs. In this paper we focus on the examples of what we call the alternating type. We classify all such examples admitting a vertex-transitive subgroup of the automorphism group of the graph preserving the corresponding $2$-factor and also determine the ones for which the $2$-factor is invariant under the full automorphism group of the graph. In this way we introduce a new infinite family of cubic vertex-transitive graphs that is a natural generalization of the well-known generalized Petersen graphs as well as of the honeycomb toroidal graphs. The family contains an infinite subfamily of arc-regular examples and an infinite family of $2$-arc-regular examples. 
\end{abstract}
\hrule

\begin{quotation}
\end{quotation}

\section{Introduction}
\label{sec:Intro}

In this paper all graphs are assumed to be finite, undirected, simple and connected, unless otherwise specified. For definitions of some terms not defined in the Introduction, see Sections~\ref{sec:prelim} and~\ref{sec:two_types}.

The main motivation for the investigations undertaken in this paper is the following problem which is a generalization of the question of Bojan Mohar to the first author of this paper (see~\cite{AlsKreKho19}):

\begin{problem}
\label{pro:main}
Classify or characterize the cubic vertex-transitive graphs admitting a partition of their edge sets into a $2$-factor $\cC$ and a $1$-factor such that some vertex-transitive subgroup of the automorphism group of the graph preserves $\cC$.
\end{problem}

Once these graphs have been determined a natural next step (which leads to the original question of Mohar) of course is to determine for which of them the $2$-factor $\cC$ is preserved by the full automorphism group of the graph. Now, a vertex-transitive subgroup $G$ of the automorphism group of a cubic vertex-transitive graph $\G$ can have one, two or three orbits in its natural action on the arc set of the graph. It is easy to see that the examples where this action has only one orbit (in which case $G$ acts arc transitively on $\G$) are the only ones that do not admit a partition of the edge set into a $2$-factor and a $1$-factor that is preserved by $G$. In this sense it is thus clear that the cubic vertex-transitive graphs that admit a partition of the edge set of the graph into a $2$-factor and a $1$-factor which is preserved by the full automorphism group of the graph are precisely those that are not arc-transitive. But what we are looking for in the context of Problem~\ref{pro:main} (and the question of Mohar) is an explicit description of such graphs.

Before proceeding a slight digression is in order. In~\cite{PotSpiVer13} a complete list of connected cubic vertex-transitive graphs of order up to $1280$ was obtained. In that paper a detailed description of the methods used to compute the graphs from each of the above mentioned three cases (depending on the number of orbits of the automorphism group on the arc set of the graph) was given. In particular, it was shown that with the exception of the girth $4$ examples (which correspond to the graphs known as prisms and M\"obius ladders) the cubic vertex-transitive graphs for which the automorphism group has two orbits on the arc set correspond naturally to certain arc-transitive quartic graphs admitting what are called arc-transitive cycle decompositions (see~\cite{MikPotWil08}). In the investigations concerning Problem~\ref{pro:main} one thus could take the approach of considering the so-called cubic graphical regular representations (the graphs whose automorphism group is regular) and the ones corresponding to these quartic arc-transitive cycle decompositions separately. There are also numerous papers on cubic vertex-transitive (and arc-transitive) graphs in the literature that one could try to use (see for instance~\cite{ConHujKutMar21, EibJajSpa19, FenKutMarYan20, PotSpiVer15, PotTol20, PotVid22} for some of the more recent ones). However, in obtaining the explicit descriptions of the graphs we are looking for this does not seem to be of considerable help, which is why we decided to take a ``direct approach''.

The first steps in the investigation of the above mentioned question of Mohar were taken in~\cite{AlsKreKho19, AlsDobKhoSpa22}, where the examples for which the corresponding $2$-factor $\cC$ consists of one or two cycles, respectively, were analyzed and the ones for which the full automorphism group preserves $\cC$ were classified. The analysis of the case where $\cC$ has two cycles was considerably more difficult than the one where $\cC$ consists of a single (Hamilton) cycle. This suggests that taking this further by considering the examples for which $\cC$ consists of a small number of cycles (say three, four or five) probably does not make much sense. On the other hand, we think it is worth considering the graphs for which we do not restrict the number of cycles in $\cC$ but instead insist that each member of $\cC$ is only adjacent to two other members of $\cC$ (in the sense that $C, C' \in \cC$ are {\em adjacent} if there exist adjacent vertices $v \in V(C)$ and $v' \in V(C')$). We say that such a graph is of {\em cycle quotient type} with respect to $\cC$. 

To explain why we think this is worth pursuing we present the data obtained by performing a computer search (using {\sc Magma}~\cite{magma}) within the above mentioned census~\cite{PotSpiVer13} up to order $500$. The meaning of the notation in the table given below is as follows. For a given integer $n$ let $\cA_n$ be the set of all connected cubic vertex-transitive but not arc-transitive graphs up to order $n$ (recall that these are precisely the graphs $\G$ that admit an $\Aut(\G)$-invariant partition of their edge sets into a $2$-factor and a $1$-factor). Then let $\cA_n^{1} \subseteq \cA_n$ be the subset of the graphs $\G$ that admit an $\Aut(\G)$-invariant partition of their edge sets into a Hamilton cycle and a $1$-factor, let $\cA_n^{2} \subseteq \cA_n$ be the subset of the graphs $\G$ that admit an $\Aut(\G)$-invariant  partition of their edge sets into a $2$-factor with two cycles and a $1$-factor, and let $\cA_n^{2*} = \cA_n^{2} \setminus \cA_n^{1}$. Next, let $\cA_n^{cq} \subseteq \cA_n$ be the subset of the graphs $\G$ that admit a partition of their edge sets into an $\Aut(\G)$-invariant $2$-factor $\cC$ and a $1$-factor such that $|\cC| \geq 3$ and $\G$ is of cycle quotient type with respect to $\cC$. Finally, let $\cA_n^{cq*} = \cA_n^{cq} \setminus (\cA_n^{1} \cup \cA_n^{2})$ and $\cA_n^* = \cA_n \setminus (\cA_n^{1} \cup \cA_n^{2} \cup \cA_n^{cq})$. The obtained data is given in the following table.  
$$
\begin{array}{c||@{\quad}c@{\quad}|@{\quad}c@{\quad}|@{\quad}c@{\quad}|@{\quad}c@{\quad}|@{\quad}c@{\quad}|@{\quad}c@{\quad}|@{\quad}c}
n & |\cA_n| & |\cA_n^{1}| & |\cA_n^{2}| & |\cA_n^{2*}| & |\cA_n^{cq}| & |\cA_n^{cq*}| & |\cA_n^*|\\ 
\hline 
100 & 627 & 278 & 204 & 110 & 255 & 110 & 129 \\			
200 & 2530 & 1164 & 770 & 360 & 1147 & 430 & 576 \\		
300 & 5657 & 2667 & 1727 & 764 & 2640 & 901 & 1325 \\	
400 & 10360 & 4785 & 3046 & 1312 & 4835 & 1614 & 2649 \\	
500 & 15703 & 7524 & 4770 & 2017 & 7538 & 2368 & 3794 	
\end{array}
$$
As one might expect, the proportion of the size of $\cA_n^*$ in $\cA_n$ seems to slowly increase as $n$ grows and the same seems to hold for the size of $\cA_n^{cq}$ in $\cA_n$. On the other hand, the proportion of the size of $\cA_n^{cq*}$ in $\cA_n$ seems to slowly decrease. Nevertheless, it still seems that $\cA_n^{cq}$ and $\cA_n^{cq*}$ are important and substantial enough subsets of $\cA_n$ to be worth investigating.

In this paper we show that the cubic vertex-transitive graphs corresponding to Problem~\ref{pro:main} which are of cycle quotient type come in two essentially different ``flavors'' - we say that some of them are of alternating cycle quotient type and some of bialternating cycle quotient type (see Section~\ref{sec:two_types} for details). It turns out that of the $7538$ graphs $\G$ from $\cA_{500}^{cq}$ around $76$ percent of them ($5733$ to be precise) admit an $\Aut(\G)$-invariant partition of their edge sets into a $2$-factor $\cC$ and a $1$-factor such that $\G$ is of alternating cycle quotient type with respect to $\cC$, while around $28$ percent of them ($2112$ to be precise) admit an $\Aut(\G)$-invariant partition of their edge sets into a $2$-factor $\cC$ and a $1$-factor such that $\G$ is of bialternating cycle quotient type with respect to $\cC$ (meaning that $307$ of them admit both kinds of cycle quotient type partitions). 

The main aim of the present paper is to thoroughly investigate the ones of alternating cycle quotient type. Our main results are as follows. We give a complete classification of connected cubic graphs $\G$ admitting a vertex-transitive subgroup $G \leq \Aut(\G)$ and a $G$-invariant partition of their edge sets into a $2$-factor $\cC$ and a $1$-factor such that $\G$ is of alternating cycle quotient type with respect to $\cC$ (see Theorem~\ref{the:altVT}). Moreover, we determine for which of these graphs the $2$-factor $\cC$ is preserved by $\Aut(\G)$ (see Theorem~\ref{the:alt_main}). It is interesting to note that the corresponding graphs $\cXa(m,n,k,\ell)$ from Theorem~\ref{the:altVT} are natural generalizations of both the well-known generalized Petersen graphs~\cite{FruGraWat70} and the honeycomb toroidal graphs~\cite{Als21}. What is more, this family of graphs contains infinite subfamilies of arc-regular graphs, $2$-arc-regular graphs, as well as graphs for which the automorphism group has two orbits on the arc set but with vertex-stabilizers of order more than $2$. In this sense the family $\cXa(m,n,k,\ell)$ is interesting on its own and may well be the subject of further investigations on cubic vertex-transitive graphs.

\section{Preliminaries}
\label{sec:prelim}

For a graph $\G$ we denote its vertex set and edge set by $V(\G)$ and $E(\G)$, respectively, where we sometimes simply write $V$ and $E$ if the graph $\G$ is clear from the context. A subgroup $G$ of the automorphism group $\Aut(\G)$ of a graph $\G$ is said to be {\em vertex-transitive}, {\em edge-transitive} or {\em $s$-arc-transitive} on $\G$, where $s \geq 1$, if the natural action of $G$ on $V(\G)$, $E(\G)$ or the set of all $s$-arcs, respectively, is transitive (where an {\em $s$-arc} is a sequence of vertices $(v_0,v_1,\ldots , v_s)$ of $\G$ for which any two consecutive vertices are adjacent and any three consecutive vertices are pairwise distinct). Moreover, if the action of $G$ on the set of $s$-arcs is regular, we say that $G$ is {\em $s$-arc-regular} on $\G$. In the case that $G = \Aut(\G)$ we say that the graph $\G$ is {\em vertex-transitive}, {\em edge-transitive}, {\em $s$-arc-transitive} or {\em $s$-arc-regular}, respectively, where we abbreviate the terms $1$-arc-transitive and $1$-arc-regular to arc-transitive and arc-regular, respectively.

Throughout the paper the Greek letters $\alpha, \beta, \gamma, \rho, \theta, \eta$ and $\tau$ will denote automorphisms of graphs while $\delta$ and $\varepsilon$ will usually be integers from $\{0,1\}$ or $\{-1,1\}$. This should cause no confusion.

We will constantly be working with elements of the residue class ring $\ZZ_n$ for some positive integer $n$, as well as with integers. Quite often we will write things such as $t\ell + 2 = 0$, where $\ell$ will be an element of $\ZZ_n$ and $t$ will be an integer. By this we mean that the corresponding equality holds in $\ZZ_n$, or in other words that viewing $\ell$ as an integer (any one in its residue class) $t\ell + 2$ is divisible by $n$. We also make the convention that for an integer (or an element of $\ZZ_n$), say $k$, we abbreviate things like $k \in \{-2,2\}$ to $k = \pm 2$ and $k \notin \{-2,2\}$ to $k \neq \pm 2$. 

A family of graphs that will play an important role in this paper is that of the honeycomb toroidal graphs $\HTG(m,n,\ell)$. Namely, this is the most natural family of examples of connected cubic graphs $\G$ admitting a partition of their edge sets into a $2$-factor $\cC$ and a $1$-factor such that $\G$ is of alternating cycle quotient type with respect to $\cC$ and some vertex-transitive subgroup of $\Aut(\G)$ preserves $\cC$. These graphs have been extensively studied (see for instance the recent survey~\cite{Als21}). For instance, it is known that they admit a regular subgroup of the automorphism group (isomorphic to a generalized dihedral group)~\cite{AlsDea09} and their automorphism groups have also been determined~\cite{Spa22}. To avoid repetition let us simply say that the {\em honeycomb toroidal graph} $\HTG(m,n,\ell)$, where $m \geq 1$, $n \geq 4$ is even and $\ell \in \ZZ_n$ is of the same parity as $m$, is precisely what we get in our Construction~\ref{cons:alt} (where we naturally extend that construction to allow for $m \in \{1,2\}$) by taking $\cXa(m,n,S,\ell)$ with $S$ consisting of $m$ copies of $\{-1,1\}$ (but see also~\cite{Als21}).

\section{The two cycle quotient types}
\label{sec:two_types}

We now describe the general setting and introduce some corresponding notation that we will be working with throughout the paper. Unless otherwise specified our graph $\G$ will always be a connected cubic vertex-transitive graph admitting a partition of its edge set into a $2$-factor $\cC$ and a $1$-factor $\cI$. So, whenever we speak of $\cC$ or $\cI$ we are referring to this chosen $2$-factor or $1$-factor of $\G$, respectively (note that a given $\G$ may admit several different partitions of its edge set into a $2$-factor and a $1$-factor). We will usually also assume that there exists some vertex-transitive subgroup $G \leq \Aut(\G)$ preserving the $2$-factor $\cC$ (which of course happens if and only if it preserves the $1$-factor $\cI$). In this case we will say that this partition is {\em $G$-invariant}.

Let $\G$, $\cC$ and $G$ be as in the previous paragraph. We let the {\em quotient graph} $\G_\cC$ of $\G$ {\em with respect to} $\cC$ be the simple graph with vertex set $\cC$ in which two distinct vertices $C$ and $C'$ are adjacent if and only if there is a vertex $v$ of $C$ that is adjacent to at least one vertex $v'$ of $C'$. Note that because of the assumption on $\cC$ and $G$ the graph $\G_\cC$ is vertex-transitive (and thus regular). Even though we will only be working with examples with $|\cC| \geq 3$ in this paper we mention that by our definition $\G_\cC$ is the complete graph $K_1$ if and only if $|\cC| = 1$ (the corresponding graphs were studied in~\cite{AlsKreKho19}), while $\G_\cC$ is the complete graph $K_2$ if and only if $|\cC| = 2$ (the corresponding graphs were studied in~\cite{AlsDobKhoSpa22}). As already stated in the Introduction, we say that $\G$ is of {\em cycle quotient type with respect to} $\cC$ if $\G_\cC$ is a cycle (or in other words if each $C \in \cC$ has precisely two neighbors in $\G_\cC$). It is these kind of graphs that we want to study in this paper. 

Suppose $\G$ is of cycle quotient type with respect to $\cC$ and suppose a vertex-transitive subgroup $G \leq \Aut(\G)$ preserves $\cC$. Since in this case $\cC$ consist of at least three cycles and $\G$ is cubic, the cycles from $\cC$ are all induced cycles of the same length, which we denote throughout the paper by $n$. Moreover, each vertex $v$ of $\G$ has precisely one neighbor which is in a different member of $\cC$ than $v$. We call this neighbor of $v$ the {\em outside neighbor} of $v$ and denote it by $\No{v}$. We let $m = |\cC|$ (and so $\G$ is of order $mn$) and we denote the members of $\cC$ by $C_i$, $i \in \ZZ_m$, where we assume that in $\G_\cC$ each $C_i$ is adjacent to $C_{i-1}$ and $C_{i+1}$, computation of indices being performed modulo $m$. Moreover, we denote the vertices of $\G$ by $\vv{i,j}$, $i \in \ZZ_m$, $j \in \ZZ_n$, in such a way that $V_i = \{\vv{i,j} \colon j \in \ZZ_n\}$ is the vertex set of $C_i$ for each $i \in \ZZ_m$. Throughout the paper all computations of the indices for $\vv{i,j}$ are thus to be performed modulo $m$ in the first and modulo $n$ in the second component. 

With no loss of generality we assume that $\vv{0,j} \sim \vv{0,j+1}$ for all $j \in \ZZ_n$ and that $\vv{0,0} \sim \vv{1,0}$. For each $i \in \ZZ_m$ we denote the setwise stabilizer of $V_i$ in $G$ by $G_i$. Since $G$ is vertex-transitive and preserves $\cC$, the restriction of the action of $G_0$ to $V_0$ is transitive, and so the corresponding group is a transitive subgroup of the dihedral group $D_n$ of order $2n$. Since $D_n$ has at most two proper transitive subgroups, namely the cyclic group of order $n$ and the dihedral group $D_{n/2}$ (when $n$ is even), it is clear that $G$ contains an automorphism $\rho$ preserving $V_0$ and whose restriction to $V_0$ is the $2$-step rotation of $C_0$, where
\begin{equation}
\label{eq:rho}
	\rho(\vv{0,j}) = \vv{0,j+2}\  \text{for\ each}\  j \in \ZZ_n.
\end{equation}  
Consider the vertex $\vv{0,0}$ of $\G$ and recall that we have decided that $\No{\vv{0,0}} = \vv{1,0} \in V_1$. There are then two essentially different possibilities. The first is that none of $\No{\vv{0,1}}$ and $\No{\vv{0,-1}}$ is in $V_1$. In this case we say that $\G$ is of {\em alternating} cycle quotient type with respect to $\cC$. In the other case we say that $\G$ is of {\em bialternating} cycle quotient type with respect to $\cC$ (note that if say $\No{\vv{0,1}} \in V_1$ then for $\G$ to be of cycle quotient type with respect to $\cC$ we require that $\No{\vv{0,2}}, \No{\vv{0,3}} \in V_{m-1}$, and so the outside neighbors of consecutive pairs of vertices of $C_0$ alternate between being in $V_1$ and $V_{m-1}$). As announced in the Introduction it is the aim of this paper to classify the graphs of alternating cycle quotient type and determine for which of them the corresponding $2$-factor $\cC$ is preserved by the whole $\Aut(\G)$. We start with the following useful observation.

\begin{proposition}
\label{pro:quo_group} 
Let $\G$ be a connected cubic graph admitting a partition of its edge set into a $2$-factor $\cC$ and a $1$-factor $\cI$ such that $|\cC| \geq 3$ and that the quotient graph $\G_\cC$ with respect to $\cC$ is a cycle. Suppose there exists a vertex-transitive subgroup $G \leq \Aut(\G)$ preserving this partition. Then the group induced by the action of $G$ on $\G_\cC$ is the dihedral group of order $2m$, where $m = |\cC|$.
\end{proposition}

\begin{proof}
Let $i \in \ZZ_m$ and let $\vv{i,j}$ be a vertex of $\G$ such that $\No{\vv{i,j}} \in V_{i+1}$. Since $G$ is vertex-transitive on $\G$ and preserves $\cC$, there exists an $\eta_i \in G$ mapping $\vv{i,j}$ to $\No{\vv{i,j}}$. Since $\vv{i,j}$ and $\No{\vv{i,j}}$ are the outside neighbors of each other, it follows that $\eta_i$ interchanges $\vv{i,j}$ and $\No{\vv{i,j}}$, and so its induced action on $\G_\cC$ interchanges $C_i$ and $C_{i+1}$. Since this holds for each $i \in \ZZ_m$ and $\G_\cC$ is a cycle of length $m$ it follows that the group induced by the action of $G$ on $\G_\cC$ is indeed the full dihedral group $D_m$. 
\end{proof}

\section{The graphs $\cXa(m,n,S,\ell)$}
\label{sec:alt_graphs}

Let $\G$, $\cC$, $\cI$, $G \leq \Aut(\G)$, $n$ and $m$ be as in the previous section where we assume that $\G$ is of alternating cycle quotient type with respect to $\cC$. Recall that $\No{\vv{0,0}} \in V_1$ and $\No{\vv{0,1}}, \No{\vv{0,-1}} \in V_{m-1}$. Let $\rho \in G$ be such that~\eqref{eq:rho} holds and note that (since $\rho(\vv{0,-1}) = \vv{0,1}$) it preserves each of $V_0$ and $V_{m-1}$. By Proposition~\ref{pro:quo_group} it then preserves each $V_i$, $i \in \ZZ_m$. Moreover, applying $\rho$ we see that $n$ must be even, say $n = 2n'$ for some $n' \geq 2$, and 
$$
	\No{\vv{0,j}} \in V_1 \iff 2 \mid j \quad \text{and} \quad \No{\vv{0,j}} \in V_{m-1} \iff 2 \nmid j.
$$
With no loss of generality we can assume that $\vv{0,j} \sim \vv{1,j}$ for all even $j \in \ZZ_n$. Therefore, $\No{\vv{1,j}} \in V_2$ for all odd $j \in \ZZ_n$, and so we can again assume that $\vv{1,j} \sim \vv{2,j}$ for all odd $j \in \ZZ_n$. Continuing in this way we can thus assume that for each $i$ with $0 \leq i \leq m-2$ we have that
\begin{equation}
\label{eq:flat_alt}
	\vv{i,j} \sim \vv{i+1,j}\ \text{for}\ \text{all}\ j \in \ZZ_n\ \text{with}\ j \equiv i \pmod{2}.
\end{equation} 
To fully describe $\G$ we now only need to describe the edges of the $m-1$ cycles $C_i$, $1 \leq i \leq m-1$, and the edges of $\cI$ connecting the vertices $\vv{0,j} \in V_0$ with $j$ odd to those of $V_{m-1}$.

Recall that none of the two neighbors of $\vv{0,0}$ in $C_0$ has its outside neighbor in the same set $V_i$ as $\vv{0,0}$. As $G$ preserves $\cC$ and is vertex-transitive it thus follows that for each $i \in \ZZ_m$ and each pair of consecutive vertices of $C_i$ their outside neighbors are in different members of $\cC$. Therefore, if for some $i \in \ZZ_m$ and $j,j' \in \ZZ_n$ the vertices $\vv{i,j}$ and $\vv{i,j'}$ are adjacent, then $j$ and $j'$ are of different parity. This proves the first half of the following result.

\begin{proposition}
\label{pro:alt_n=4}
Let $\G$ be a connected cubic vertex-transitive graph admitting a partition of its edge set into a $2$-factor $\cC$ and a $1$-factor $\cI$ such that $\G$ is of alternating cycle quotient type with respect to $\cC$. If $\cC$ is preserved by a vertex-transitive subgroup $G \leq \Aut(\G)$ and the members of $\cC$ are $4$-cycles, then $\G$ is isomorphic to the honeycomb toroidal graph $\mathrm{HTG}(m,4,\ell)$, where $m = |\cC|$ and $\ell \in \{0,1\}$ is of the same parity as $m$. Moreover, $\cC$ is preserved by $\Aut(\G)$.
\end{proposition}

\begin{proof}
That $\G$ is isomorphic to $\mathrm{HTG}(m,4,\ell)$ for the unique $\ell \in \{0,1\}$ which is of the same parity as $m$ follows from the above discussion and the fact that $\mathrm{HTG}(m,4,\ell) \cong \mathrm{HTG}(m,4,\ell + 2)$ (we can simply exchange the roles of $\vv{0,1}$ and $\vv{0,3}$). That the partition is $\Aut(\G)$-invariant follows from the fact that the cycles from the $2$-factor $\cC$ are $4$-cycles while the edges of the $1$-factor do not lie on $4$-cycles (as $|\cC| \geq 3$).
\end{proof}

In view of the above proposition we can assume henceforth in this section that $n' \geq 3$ (recall that $n = 2n'$). Let $\rho \in G$ be such that its action on $V_0$ is as in~\eqref{eq:rho}. Since it preserves the $1$-factor $\cI$, it maps each $\vv{1,j}$ with $j$ even to $\vv{1,j+2}$. Thus $n' \geq 3$ implies that the restriction of the action of $\rho$ to $C_1$ is a rotation having two orbits of size $n'$ on $V_1$. With no loss of generality we can thus assume that the vertices $\vv{1,j}$ with $j$ odd have been labeled in such a way that $\rho(\vv{1,j}) = \vv{1,j+2}$ for all $j \in \ZZ_n$. We now see that $\rho(\vv{2,j}) = \vv{2,j+2}$ for all odd $j \in \ZZ_n$, and so we can assume that the vertices $\vv{2,j}$ with $j$ even have been labeled in such a way that $\rho(\vv{2,j}) = \vv{2,j+2}$ for all $j \in \ZZ_n$. Continuing in this way we thus finally see that we can assume that the vertices of $\G$ have been labeled in such a way that
\begin{equation}
\label{eq:rho_alt}
	\rho(\vv{i,j}) = \vv{i,j+2} \ \text{for}\ \text{all}\ i \in \ZZ_m\ \text{and}\ j \in \ZZ_n.
\end{equation}
It now also follows that there exists a unique $\ell \in \ZZ_n$ such that 
\begin{equation}
\label{eq:ell_alt}
	\vv{m-1,j} \sim \vv{0,j+\ell}\ \text{for}\ \text{all}\ j \in \ZZ_n\ \text{with}\ j \equiv m-1 \pmod{2}.
\end{equation}
Observe that $\ell$ must be of the same parity as $m$ so that $j+\ell$ will be odd and thus $\vv{0,j+\ell}$ will indeed be a vertex having its outside neighbor in $V_{m-1}$. Finally, the discussion in the paragraph preceding Proposition~\ref{pro:alt_n=4} together with~\eqref{eq:rho_alt} implies that for each $i \in \ZZ_m$ there exist distinct odd $a_i,b_i \in \ZZ_n$ such that 
\begin{equation}
\label{eq:vert_alt}
	\vv{i,j} \sim \vv{i,j+a_i}, \vv{i,j+b_i} \ \text{for}\ \text{all}\ j \in \ZZ_n\ \text{with}\ j \equiv i \pmod{2}.
\end{equation}
Of course, $\{a_0, b_0\} = \{-1,1\}$. For the edges connecting the vertices of $V_i$ to form an $n$-cycle we of course require $\gcd(b_i-a_i,n) = 2$ to hold. We now have a complete description of the graph $\G$.

\begin{construction}
\label{cons:alt}
Let $m, n$ be integers where $n \geq 4$ is even and $m \geq 3$ and let $\ell \in \ZZ_n$ be of the same parity as $m$. Furthermore, set $a_0 = 1$ and $b_0 = -1$ and for each $i$ with $1 \leq i \leq m-1$ let $a_i,b_i \in \ZZ_n$ be distinct odd elements such that $\gcd(b_i-a_i,n)=2$. Then the graph $\cXa(m,n,S,\ell)$, where $S = [\{a_0,b_0\}, \{a_1,b_1\},\ldots , \{a_{m-1},b_{m-1}\}]$, is the graph of order $mn$ with vertex set consisting of all $\vv{i,j}$, $i \in \ZZ_m$, $j \in \ZZ_n$, and adjacencies given in~\eqref{eq:flat_alt}, \eqref{eq:ell_alt} and~\eqref{eq:vert_alt}.
\end{construction}

\begin{proposition}
\label{pro:theQgraphs}
Let $\G$ be a connected cubic vertex-transitive graph admitting a partition of its edge set into a $2$-factor $\cC$ and a $1$-factor such that $\G$ is of alternating cycle quotient type with respect to $\cC$. If $\cC$ is preserved by some vertex-transitive subgroup of $\Aut(\G)$ then $\G$ is isomorphic to a graph $\cXa(m,n,S,\ell)$ from Construction~\ref{cons:alt}, where $m = |\cC|$. 
\end{proposition}

The sequence $S$ from the above construction is called the {\em signature} of the graph $\cXa(m,n,S,\ell)$. In the next section we show that it suffices to consider only signatures of a very particular kind (see Corollary~\ref{cor:alt_iso} and Proposition~\ref{pro:alt_signatures}). We also make the agreement that whenever we have a graph $\G = \cXa(m,n,S,\ell)$ from Construction~\ref{cons:alt} we let $V_i = \{\vv{i,j} \colon j \in \ZZ_n\}$ and $C_i$ be the cycle induced on $V_i$ for each $i \in \ZZ_m$. Moreover, we set $\cC = \{C_i \colon i \in \ZZ_m\}$ and we let $\cI$ be the set of all edges arising from adjacencies in~\eqref{eq:flat_alt} and~\eqref{eq:ell_alt}. We call the edges from $\cI$ the {\em links} of $\G$ and the remaining ones (that is, those from the cycles in $\cC$) the {\em non-links} of $\G$.

\section{Simplifying the signatures}
\label{sec:alt_signatures}

We first record some useful isomorphisms between graphs from Construction~\ref{cons:alt}.

\begin{lemma}
\label{le:alt_iso}
Let $\G = \cXa(m,n,S,\ell)$, where the parameters satisfy all the assumptions from Construction~\ref{cons:alt}. Then for each even $q \in \ZZ_n$ both of the following hold:
\begin{itemize}
\item[(i)] For any $i$ with $1 \leq i \leq m-2$ we have that $\G \cong \cXa(m,n,S',\ell)$, where the signature $S'$ is obtained from $S$ by replacing $\{a_i,b_i\}$ and $\{a_{i+1},b_{i+1}\}$ by $\{a_{i}-q, b_{i}-q\}$ and $\{a_{i+1}+q, b_{i+1}+q\}$, respectively.
\item[(ii)] $\G \cong \cXa(m,n,S'',\ell-q)$, where the signature $S''$ is obtained from $S$ by replacing $\{a_{m-1},b_{m-1}\}$ by $\{a_{m-1}-q, b_{m-1}-q\}$.
\end{itemize}
\end{lemma}

\begin{proof}
To confirm the claim from item (i) simply relabel each $u_{i,j}$ and $u_{i+1,j}$ where $j \equiv i \pmod{2}$ by $u_{i,j+q}$ and $u_{i+1,j+q}$, respectively. To confirm the claim from item (ii) simply relabel the vertices $u_{m-1,j}$ with $j \equiv m-1 \pmod{n}$ by $u_{m-1,j+q}$. We leave the details to the reader.
\end{proof}

\begin{corollary}
\label{cor:alt_iso}
Let $\G = \cXa(m,n,S,\ell)$, where the parameters satisfy all the assumptions from Construction~\ref{cons:alt}. Then $\G \cong \cXa(m,n,S',\ell')$ for some $\ell' \in \ZZ_n$ and a signature $S'$ such that for each $i \in \ZZ_m$ we have that $a'_i+b'_i = 0$.
\end{corollary}

\begin{proof}
Let $i$, $1 \leq i \leq m-1$, be the smallest integer such that $a_i+b_i \neq 0$ (if there is no such $i$ there is nothing to prove). We claim that there exists an even $q \in \ZZ_n$ such that $a_i-q+b_i-q = 0$, or equivalently that $2q = a_i + b_i$. Note that if such an even $q$ does indeed exist, then applying Lemma~\ref{le:alt_iso} (the first or the second item, depending on whether $i < m-1$ or not, respectively) will result in a graph $\cXa(m,n,S',\ell')$ for which the smallest $i'$ with $a_{i'}+b_{i'} \neq 0$ (if it exists at all) is larger than $i$. Doing this at most $m-1$ times we thus eventually arrive at a signature of the desired form. Of course, in the last step the parameter $\ell$ may change.

To prove our claim about the existence of an appropriate even $q \in \ZZ_n$ observe first that as $a_i$ and $b_i$ are both odd and $n$ is even, $a_i+b_i$ is even. There thus exists a $q \in \ZZ_n$ such that $2q = a_i+b_i$. If $q$ is even, we are done. If it is not, then the fact that $\gcd(b_i-a_i,n) = 2$ implies that $n/2$ is odd, and so we can simply replace $q$ by $q + n/2$. 
\end{proof}

By Corollary~\ref{cor:alt_iso} we can restrict our attention to graphs $\G = \cXa(m,n,S,\ell)$, where the signature $S$ is of the form $S = [\{\pm 1\}, \{\pm k_1\}, \ldots , \{\pm k_{m-1}\}]$ for some odd $k_1, k_2, \ldots , k_{m-1} \in \ZZ_n$, where for each $i$ we have that $\gcd(2k_i, n) = 2$, which is of course equivalent to $\gcd(k_i,n) = 1$. Observe that this implies that $\vv{i,j} \sim \vv{i,j\pm k_i}$ for all $i \in \ZZ_m$ and $j \in \ZZ_n$ (where we let $k_0 = 1$). Whenever the signature is in such a form we simply write it as $[1,k_1,k_2,\ldots , k_{m-1}]$ (where we can in fact assume that $k_i < n'$ for all $i \in \ZZ_m$). For ease of reference we now record our assumption that we will be working with in some of the next results.

\begin{assumption}
\label{assump_alt}
We assume that $\G = \cXa(m,n,S,\ell)$ for some $m \geq 3$, an even $n \geq 6$ and $\ell \in \ZZ_n$ of the same parity as $m$ and where the signature $S$ is of the form $S = [k_0,k_1,\ldots , k_{m-1}]$ for $k_0 = 1$ and some $k_1, k_2, \ldots , k_{m-1} \in \ZZ_n$ which are all coprime to $n$. We also assume that there exists a vertex-transitive subgroup $G \leq \Aut(\G)$ preserving the 2-factor $\cC$ consisting of the $m$ cycles of length $n$ induced on the sets $V_i$, $i \in \ZZ_m$. 
\end{assumption}

We now investigate in what way the assumption on the existence of the group $G$ restricts the signature $S$ and the parameter $\ell$. We achieve this via elements of $G$ guaranteed by Proposition~\ref{pro:quo_group}.

\begin{lemma}
\label{le:alt_regular}
Let $\G = \cXa(m,n,S,\ell)$ be as in Assumption~\ref{assump_alt}. Then there exists a nontrivial automorphism of $\G$ preserving the $2$-factor $\cC$ and fixing $\vv{0,0}$ if and only if $2\ell = 0$. 
\end{lemma}

\begin{proof}
Suppose $\eta \in \Aut(\G)$ is a nontrivial automorphism preserving $\cC$ and fixing $\vv{0,0}$. It then either also fixes $\vv{0,1}$ or interchanges it with $\vv{0,-1}$. In the former case it must fix the whole $C_0$ pointwise and then the fact that $n \geq 6$ implies that $\eta$ fixes at least three vertices of the cycle $C_1$ (recall that~\eqref{eq:flat_alt} holds), and thus has to fix it pointwise. Continuing in this way we thus see that $\eta$ is the identity, a contradiction. 

It thus follows that $\eta$ interchanges $\vv{0,1}$ with $\vv{0,-1}$, and consequently also interchanges $\vv{0,j}$ with $\vv{0,-j}$ for each $j \in \ZZ_n$. Considering outside neighbors it thus follows that $\eta$ interchanges $\vv{1,j}$ with $\vv{1,-j}$ for each even $j \in \ZZ_n$. Let $j \in \ZZ_n$ be odd. Since $n \geq 6$ and $C_1$ is an $n$-cycle, the vertex $\vv{1,j}$ is the unique common neighbor of $\vv{1,j-k_1}$ and $\vv{1,j+k_1}$, which are mapped by $\eta$ to $\vv{1,-j+k_1}$ and $\vv{1,-j-k_1}$, respectively (recall that $k_1$ is odd, and so $j-k_1$ and $j+k_1$ are even). Since the unique common neighbor of these two vertices is $\vv{1,-j}$,  it follows that $\eta$ interchanges $\vv{1,j}$ with $\vv{1,-j}$ for each $j \in \ZZ_n$. Continuing in this way we finally see that $\eta$ also interchanges $\vv{m-1,j}$ with $\vv{m-1,-j}$ for each $j \in \ZZ_n$. For $j \in \ZZ_n$ of the same parity as $m-1$ we thus have that $\eta$ maps the pair of adjacent vertices $\vv{m-1,j}$ and $\vv{0,j+\ell}$ to the pair $\vv{m-1,-j}$ and $\vv{0,-j-\ell}$, which are adjacent if and only if $-j+\ell = -j-\ell$. This clearly holds if and only if $2\ell = 0$. 

To prove the converse one simply needs to show that if $2\ell = 0$ then the permutation $\eta$ interchanging $\vv{i,j}$ with $\vv{i,-j}$ for each $i \in \ZZ_m$ and each $j \in \ZZ_n$ is an automorphism of $\G$, which is easy to do and is left to the reader.
\end{proof}

\begin{lemma}
\label{le:alt_swap}
Let $\G = \cXa(m,n,S,\ell)$ be as in Assumption~\ref{assump_alt}. Then $2k_i^2 = \pm 2$ holds for all $i \in \ZZ_m$. Moreover, if $2\ell \neq 0$, then $2k_i^2 = 2$ holds for all $i \in \ZZ_m$.
\end{lemma}

\begin{proof}
Let $i \in \ZZ_m$, where $i \neq m-1$, and let $\delta$ be $0$ or $1$ depending on whether $i$ is even or odd, respectively. Note that by~\eqref{eq:flat_alt} this implies that $\No{\vv{i,\delta}} = \vv{i+1,\delta}$. Since $G$ is vertex-transitive, there thus exists an $\alpha \in G$ interchanging $\vv{i,\delta}$ with $\vv{i+1,\delta}$. Hence $\alpha^2$ fixes $\vv{i,\delta}$, and so Lemma~\ref{le:alt_regular} implies that $\alpha^2$ either fixes each vertex of $\G$ or $2\ell = 0$ in which case a similar argument as in the proof of Lemma~\ref{le:alt_regular} shows that $\alpha^2$ interchanges $\vv{i',j}$ with $\vv{i',-j+2\delta}$ for each $i' \in \ZZ_m$ and $j \in \ZZ_n$. Let $\varepsilon \in \{-1,1\}$ be such that $\alpha(\vv{i,\delta + k_i}) = \vv{i+1,\delta + \varepsilon k_{i+1}}$. Then 
$$
	\alpha(\vv{i,\delta + jk_i}) = \vv{i+1,\delta + \varepsilon jk_{i+1}}\ \text{for}\ \text{all}\ j \in \ZZ_n.
$$
By~\eqref{eq:flat_alt} we thus find that
$$
	\alpha^2(\vv{i, \delta + 2k_i^2}) = \alpha(\vv{i+1, \delta + \varepsilon 2k_ik_{i+1}}) = \vv{i, \delta + 2k_{i+1}^2}, 
$$
and so either $2k_{i+1}^2 = 2k_{i}^2$ or $2\ell = 0$ and $2k_{i+1}^2 = -2k_i^2$. Since $k_0 = 1$, we are done.
\end{proof}

\begin{proposition}
\label{pro:alt_gamma}
Let $\G = \cXa(m,n,S,\ell)$ be as in Assumption~\ref{assump_alt}. Then precisely one of the following holds:
\begin{itemize}
\itemsep = 0pt
\item $2k_i = \pm 2$ holds for all $i \in \ZZ_m$, or
\item $m$ is even, $2k_1 \neq \pm 2$, and $2k_i = \pm 2$ for each even $i \in \ZZ_m$, while $2k_i = \pm 2k_1$ for each odd $i \in \ZZ_m$.
\end{itemize}
\end{proposition}

\begin{proof}
Let $i \in \ZZ_m$ be such that $i \leq m-3$ and let $\delta$ be $0$ or $1$, depending on whether $i$ is even or odd, respectively. Since $G$ is vertex-transitive there exists a $\gamma \in G$ mapping $\vv{i,\delta}$ to $\vv{i+1,\delta + 1}$. Since $\No{\vv{i,\delta}} = \vv{i+1,\delta}$ and $\No{\vv{i+1,\delta + 1}} = \vv{i+2,\delta + 1}$, we have that $\gamma(\vv{i+1,\delta}) = \vv{i+2,\delta + 1}$. The induced action of $\gamma$ on the quotient graph $\G_\cC$ is thus a $1$-step rotation mapping $V_{i'}$ to $V_{i'+1}$ for each $i' \in \ZZ_m$. Let $\varepsilon \in \{-1,1\}$ be the unique element such that $\gamma(\vv{i,\delta + k_i}) = \vv{i+1,\delta+1+\varepsilon k_{i+1}}$. Then Lemma~\ref{le:alt_swap} yields 
$$
\gamma(\vv{i,\delta + 2k_{i+1}k_i}) = \vv{i+1,\delta+1+2\varepsilon k_{i+1}^2} \in \{\vv{i+1,\delta+3}, \vv{i+1,\delta-1}\}.
$$
Applying~\eqref{eq:flat_alt} we thus find that $\gamma(\vv{i+1,\delta+2k_ik_{i+1}}) \in \{\vv{i+2,\delta+3}, \vv{i+2,\delta-1}\}$. On the other hand, since $\gamma(\vv{i+1,\delta}) = \vv{i+2,\delta+1}$ and there is a natural walk of length $2k_i$ in $C_{i+1}$ from $\vv{i+1,\delta}$ to $\vv{i+1,\delta+2k_ik_{i+1}}$ we also have that 
$$
\gamma(\vv{i+1,\delta+2k_ik_{i+1}}) \in \{\vv{i+2,\delta+1+2k_ik_{i+2}}, \vv{i+2,\delta+1-2k_ik_{i+2}}\}.
$$
It thus follows that $2k_ik_{i+2} = \pm 2$. Multiplying by $k_i$ and applying Lemma~\ref{le:alt_swap} we thus find that $2k_{i+2} = \pm 2k_i$ for all $i \in \ZZ_m$ with $i \leq m-3$. Consequently, $2k_i = \pm 2$ for all even $i \in \ZZ_m$, while $2k_i = \pm 2k_1$ for all odd $i \in \ZZ_m$. A completely analogous argument where we take $i = m-2$ shows that $2k_{m-2} = 2k_{m-2}k_0 = \pm 2$, and so if $m$ is odd we in fact have that $2k_i = \pm 2$ for all $i \in \ZZ_m$.
\end{proof}

We are now ready to make the final step in simplifying the signature of our graphs.

\begin{proposition}
\label{pro:alt_signatures}
Let $\G = \cXa(m,n,S,\ell)$ be as in Assumption~\ref{assump_alt}. Then we can relabel the vertices of each of the cycles of the $2$-factor $\cC$ in such a way that one of the following holds:
\begin{itemize}
\itemsep = 0pt
\item there exists some $\ell' \in \ZZ_n$ such that $\G = \HTG(m,n,\ell')$, or
\item $m$ is even and there exist $\ell' \in \ZZ_n$ and an odd $k \in \ZZ_n$ with $2k \neq \pm 2$ and $2k^2 = \pm 2$ such that $\G = \cXa(m,n,[1, k, 1, k, \ldots , 1, k],\ell')$.
\end{itemize}
\end{proposition}

\begin{proof}
We consider each of the two cases given by Proposition~\ref{pro:alt_gamma} separately. 

Suppose first that $2k_i = \pm 2$ for all $i \in \ZZ_m$. Then for each $i \in \ZZ_m$ we have that either $k_i = \pm 1$, in which case we can simply assume that $k_i = 1$, or $k_i = n' \pm 1$ (recall that $n = 2n'$). If the former holds for all $i \in \ZZ_m$, there is nothing to prove as then $\G$ is a honeycomb toroidal graph. Assume thus that this is not the case and let $i$, $1 \leq i \leq m-1$, be the smallest integer such that $k_i = n' \pm 1$. With no loss of generality we can in fact assume that $k_i = n'+1$. Since $k_i$ is odd, $n'$ is even (and thus $n$ is divisible by $4$). We can now apply Lemma~\ref{le:alt_iso} with $q = n'$ to obtain a new signature in which all the $k_{i'}$ with $i' < i$ remain $1$ and the new $k_i$ becomes $1$ as well. Clearly, the new signature still satisfies the assumption that $2k_{i'} = \pm 2$ for all $i' \in \ZZ_m$. We proceed in this way until finally reaching the situation in which all the $k_i$ are $1$, proving that $\G$ is a honeycomb toroidal graph (note that at the last step $\ell$ may change).

Suppose now that $2k_i \neq \pm 2$ for at least one $i \in \ZZ_m$. Then Proposition~\ref{pro:alt_gamma} implies that $m$ is even, $2k_1 \neq \pm 2$ and for each $i \in \ZZ_m$ either $2k_i = \pm 2$ or $2k_i = \pm 2k_1$, depending on whether $i$ is even or odd, respectively. The argument is now very similar as in the previous case. If $k_i \neq \pm 1$ for some even $i \in \ZZ_m$ or $k_i \neq \pm k_1$ for some odd $i \in \ZZ_m$, then the fact that $k_i \pm 1$ and $k_i \pm k_1$ are both even implies that $n'$ is even and that we can assume $k_i = n'+1$ or $k_i = n'+k_1$. Using Lemma~\ref{le:alt_iso} with $q = n'$ again shows that we can change the signature in such a way that $k_i$ becomes $1$ or $k_1$, depending on whether $i$ is even or odd, respectively. 
\end{proof}

Since the automorphism groups of honeycomb toroidal graphs are known~\cite{Spa22}, Proposition~\ref{pro:alt_signatures} implies that we only need to investigate the graphs $\cXa(m,n,S,\ell)$ with $m$ even and $S$ of the form $[1,k,1, \ldots , k]$, where $k \in \ZZ_n$ is such an odd element that $2k \neq \pm 2$ and $2k^2 = \pm 2$. We record the adjusted assumption for ease of reference and introduce a simplified notation for such graphs (we point out that here we do not assume that a vertex-transitive subgroup of the automorphism group preserving the $2$-factor $\cC$ exists).

\begin{assumption}
\label{assump_alt2}
We assume $\G = \cXa(m,n,S,\ell)$, where $m \geq 4$, $n \geq 6$ and $\ell \in \ZZ_n$ are all even and $S = [1,k,1,k, \ldots ,1,k]$ for an odd $k \in \ZZ_n$ such that $2k \neq \pm 2$ and $2k^2 = \pm 2$. We simplify the notation for this graph to $\cXa(m,n,k,\ell)$.
\end{assumption}

Note that for $2k \neq \pm 2$ to be possible, we require that $n \geq 10$. The smallest possible graphs satisfying Assumption~\ref{assump_alt2} are thus of order $40$ (having $m = 4$ and $n = 10$). The examples $\cXa(4,10,3,0)$ and $\cXa(4,10,3,2)$ are presented in Figure~\ref{fig:alt_example} ($C_0$ corresponds to the ``outer'' $10$-cycle). It turns out that the first of these two graphs is $2$-arc-regular, while the second is not even vertex-transitive.

\begin{figure}[h]
\begin{center}
\subfigure
{
\begin{tikzpicture}[scale = .4]
\foreach \i in {0,1,2,3}{
\pgfmathtruncatemacro{\rad}{\i+1};
\foreach \j in {0,1,2,3,4,5,6,7,8,9}{
\node[vtx, ,fill = black, inner sep = 2pt,] (A\i\j) at (360*\j/10:1.7*\rad) {};
}}
\begin{scope}[on background layer]
\foreach \j in {0,1,2,3,4,5,6,7,8,9}{
\foreach \i in {1,3}{
\draw[thick] let \n1 = {int(mod(\j+1, 10))} in (A\i\j) -- (A\i\n1);
}
\foreach \i in {0,2}{
\draw[thick] let \n1 = {int(mod(\j+3, 10))} in (A\i\j) to[bend left = 10] (A\i\n1);
}}
\foreach \j in {0,2,4,6,8}{
\foreach \i in {1,3}{
\draw[thick] let \n1 = {int(mod(\i-1, 4))} in (A\n1\j) -- (A\i\j);
}}
\foreach \j in {1,3,5,7,9}{
\draw[thick] (A2\j) -- (A1\j);
\draw[thick] (A3\j) to[bend left = 20] (A0\j);
}
\end{scope}
\end{tikzpicture}
}
\quad \quad
\subfigure
{
\begin{tikzpicture}[scale = .4]
\foreach \i in {0,1,2,3}{
\pgfmathtruncatemacro{\rad}{\i+1};
\foreach \j in {0,1,2,3,4,5,6,7,8,9}{
\node[vtx, ,fill = black, inner sep = 2pt,] (A\i\j) at (360*\j/10:1.7*\rad) {};
}}
\begin{scope}[on background layer]
\foreach \j in {0,1,2,3,4,5,6,7,8,9}{
\foreach \i in {1,3}{
\draw[thick] let \n1 = {int(mod(\j+1, 10))} in (A\i\j) -- (A\i\n1);
}
\foreach \i in {0,2}{
\draw[thick] let \n1 = {int(mod(\j+3, 10))} in (A\i\j) to[bend left = 10] (A\i\n1);
}}
\foreach \j in {0,2,4,6,8}{
\foreach \i in {1,3}{
\draw[thick] let \n1 = {int(mod(\i-1, 4))} in (A\n1\j) -- (A\i\j);
}}
\foreach \j in {1,3,5,7,9}{
\draw[thick] (A2\j) -- (A1\j);
\draw[thick] let \n1 = {int(mod(\j+2, 10))} in (A0\n1) to[bend left = 10] (A3\j);
}
\end{scope}
\end{tikzpicture}
}
\caption{The graphs $\cXa(4,10,3,0)$ and $\cXa(4,10,3,2)$.}\label{fig:alt_example}
\end{center}
\end{figure}

\section{Achieving vertex transitivity}
\label{sec:alt_VT}

In this section we determine which of the graphs $\G$ from Assumption~\ref{assump_alt2} admit a vertex-transitive subgroup of $\Aut(\G)$ preserving the $2$-factor $\cC$ thereby solving Problem~\ref{pro:main} for the case of the graphs of alternating cycle quotient type. 

\begin{proposition}
\label{pro:alt_gamma2}
Let $\G = \cXa(m,n,k,l)$ be as in Assumption~\ref{assump_alt2}. Then there exists some $\gamma \in \Aut(\G)$ preserving the $2$-factor $\cC$ and mapping $\vv{0,0}$ to $\vv{1,1}$ if and only if one of the following holds:
\begin{itemize}
\itemsep = 0pt
\item[(i)] $\ell k = \pm \ell$ and at least one of $k^2 = \pm 1$ and $4 \mid m$ holds.
\item[(ii)] $k^2 \neq \pm 1$, $m \equiv 2 \pmod{4}$ and $\ell k = n' \pm \ell$, where $n = 2n'$.
\end{itemize}
\end{proposition}

\begin{proof}
We distinguish two cases depending on whether $k^2 = \pm 1$ or not.
\smallskip

\noindent
{\sc Case 1:} $k^2 = \pm 1$.\\
Let $\delta \in \{-1,1\}$ be such that $k^2 = \delta$ (and so $\delta k^2 = 1$). Suppose there exists a $\gamma \in \Aut(\G)$ preserving $\cC$ and mapping $\vv{0,0}$ to $\vv{1,1}$. Letting $\varepsilon \in \{-1,1\}$ be such that $\gamma(\vv{0,1}) = \vv{1,1+\varepsilon k}$ we see that
\begin{equation}
\label{eq:alt_gamma_0}
	\gamma(\vv{0,j}) = \vv{1,1+\varepsilon jk}\ \text{for}\ \text{all}\ j \in \ZZ_n.
\end{equation}
Since $\gamma$ preserves the set of all links of $\G$,~\eqref{eq:flat_alt} implies that $\gamma(\vv{1,j}) = \vv{2,1+\varepsilon jk}$ for all even $j \in \ZZ_n$. In particular, $\gamma(\vv{1,0}) = \vv{2,1}$ and $\gamma(\vv{1,2k}) = \vv{2,1+2\varepsilon\delta}$, and so the common neighbor $\vv{1,k}$ of $\vv{1,0}$ and $\vv{1,2k}$ is mapped by $\gamma$ to $\vv{2,1+\varepsilon\delta}$. It now clearly follows that $\gamma(\vv{1,jk}) = \vv{2,1+\varepsilon\delta j}$ holds for all $j \in \ZZ_n$. Taking into account that $j = (\delta jk)k$ we thus see that $\gamma(\vv{1,j}) = \vv{2,1+\varepsilon jk}$ for all $j \in \ZZ_n$. Continuing in this way we find that
\begin{equation}
\label{eq:alt_gamma_1}
	\gamma(\vv{i,j}) = \vv{i+1,1+\varepsilon jk}\ \text{for}\ \text{all}\ i \in \ZZ_m \setminus \{m-1\}\ \text{and}\ j \in \ZZ_n,
\end{equation}
and  
\begin{equation}
\label{eq:alt_gamma_2}
	\gamma(\vv{m-1,j}) = \vv{0,1+\ell+\varepsilon jk}\ \text{for}\ \text{all}\ j \in \ZZ_n.
\end{equation}
Since $m$ is even, \eqref{eq:ell_alt} implies that $\vv{m-1,j} \sim \vv{0,j+\ell}$ for each odd $j \in \ZZ_n$. Since these two vertices are mapped by $\gamma$ to $\vv{0,1+\ell+\varepsilon jk}$ and $\vv{1,1+\varepsilon jk + \varepsilon \ell k}$, it thus follows that $\ell k = \varepsilon \ell$ holds. To prove the converse we simply show that if $\ell k = \varepsilon \ell$ for some $\varepsilon \in \{-1,1\}$, then the permutation $\gamma$ defined on $\G$ by~\eqref{eq:alt_gamma_1} and~\eqref{eq:alt_gamma_2} is indeed an automorphism of $\G$. We leave this easy verification to the reader.
\smallskip

\noindent
{\sc Case 2:} $k^2 \neq \pm 1$.\\
Since $2k^2 = \pm 2$ we thus have that $k^2 = n' + \delta$ for a $\delta \in \{-1,1\}$. Since $k$ is odd this implies that $n'$ is even and that $1 = \delta (k+n')k$. Suppose again that there exists a $\gamma \in \Aut(\G)$ preserving $\cC$ and mapping $\vv{0,0}$ to $\vv{1,1}$. We proceed very similarly as above, first letting $\varepsilon \in \{-1,1\}$ be such that $\gamma(\vv{0,1}) = \vv{1,1+\varepsilon k}$ and obtaining that~\eqref{eq:alt_gamma_0} holds and consequently that $\gamma(\vv{1,jk}) = \vv{2,1+\varepsilon\delta j}$ for all $j \in \ZZ_n$. Therefore, 
$$
	\gamma(\vv{1,j}) = \gamma(\vv{1,\delta j (k+n_0)k}) = \vv{2,1+\varepsilon jk + jn'}\ \text{for}\ \text{all}\ j \in \ZZ_n.
$$
By~\eqref{eq:flat_alt} we then obtain that $\gamma(\vv{2,j}) = \vv{3,1+\varepsilon jk + n'}$ holds for all odd $j \in \ZZ_n$, and consequently also 
$$
	\gamma(\vv{2,j}) = \vv{3,1+\varepsilon jk + n'}\ \text{for}\ \text{all}\ j \in \ZZ_n.
$$
Continuing in this way we finally find that
\begin{equation}
\label{eq:alt_gamma_3}
	\gamma(\vv{i,j}) = \left\{\begin{array}{ccl}
		\vv{i+1,1+\varepsilon jk} & : & i \equiv 0 \pmod{4}\\
		\vv{i+1,1+\varepsilon jk + jn'} & : & i \equiv 1 \pmod{4}\\
		\vv{i+1,1+\varepsilon jk + n'} & : & i \equiv 2 \pmod{4}\\
		\vv{i+1,1+\varepsilon jk + (j+1)n'} & : & i \equiv 3 \pmod{4},\end{array}\right.\quad i \in \ZZ_m \setminus \{m-1\},\ j \in \ZZ_n.
\end{equation}
Moreover, if $4\mid m$, then 
\begin{equation}
\label{eq:alt_gamma_4}
	\gamma(\vv{m-1,j}) = \vv{0,1+\ell+\varepsilon jk + (j+1)n'}\ \text{for}\ \text{all}\ j \in \ZZ_n,
\end{equation}
while if $m \equiv 2 \pmod{4}$, then
\begin{equation}
\label{eq:alt_gamma_5}
	\gamma(\vv{m-1,j}) = \vv{0,1+\ell+\varepsilon jk + jn'}\ \text{for}\ \text{all}\ j \in \ZZ_n.
\end{equation}
We now again consider the $\gamma$-image of a link of the form $\vv{m-1,j}\vv{0,j+\ell}$, $j \in \ZZ_n$ odd. In the case of $4 \mid m$ we find that $\ell k = \varepsilon \ell$ must hold, while in the case of $m \equiv 2 \pmod{4}$ we find that $\ell k = \varepsilon \ell + n'$, as claimed. Again, one can easily verify that the converse holds in the sense that if $4 \mid m$ and $\ell k = \varepsilon \ell$, then the mapping $\gamma$ defined by~\eqref{eq:alt_gamma_3} and~\eqref{eq:alt_gamma_4} is an automorphism of $\G$, while if $m \equiv 2 \pmod{4}$ and $\ell k = \varepsilon \ell + n'$, then the $\gamma$ defined by~\eqref{eq:alt_gamma_3} and~\eqref{eq:alt_gamma_5} is an automorphism of $\G$.
\end{proof}

We can now classify the graphs of alternating cycle quotient type that correspond to Problem~\ref{pro:main}.

\begin{theorem}
\label{the:altVT}
Let $\G$ be a connected cubic vertex-transitive graph admitting a partition of its edge set into a $2$-factor $\cC$ and a $1$-factor $\cI$ such that $\G$ is of alternating cycle quotient type with respect to $\cC$. Let $m = |\cC|$ and let $n$ be such that $\G$ is of order $mn$. Then there exists a vertex-transitive subgroup of $\Aut(\G)$ preserving $\cC$ if and only if one of the following holds:
\begin{itemize}
\itemsep = 0pt
	\item[(i)] $\G \cong \HTG(m,n,\ell)$ for some $\ell \in \ZZ_n$ of the same parity as $m$, or
	\item[(ii)] $m$ is even and $\G \cong \cXa(m,n,k,\ell)$ for some even $\ell \in \ZZ_n$ and odd $k \in \ZZ_n$ with $2k \neq \pm 2$ and $2k^2 = \pm 2$ and one of the following holds:
	\begin{itemize}
	\itemsep = 0pt
		\item[$\bullet$] $\ell k = \pm \ell$ and at least one of $k^2 = \pm 1$ and $4 \mid m$ holds, or
		\item[$\bullet$] $\ell k = n' \pm \ell$, $m \equiv 2 \pmod{4}$ and $k^2 = n' \pm 1$, where $n = 2n'$.
	\end{itemize}
\end{itemize}
\end{theorem}

\begin{proof}
The forward implication follows from Proposition~\ref{pro:alt_n=4}, Proposition~\ref{pro:theQgraphs}, Corollary~\ref{cor:alt_iso}, Proposition~\ref{pro:alt_signatures} and Proposition~\ref{pro:alt_gamma2}. For the converse we first note that the $\HTG(m,n,\ell)$ graphs indeed do admit a vertex-transitive subgroup of their automorphism group preserving $\cC$. As we already mentioned they admit a regular generalized dihedral group as a subgroup of the automorphism group and the links of these graphs correspond to a certain involution of this group. 

To complete the proof we thus only need to consider the graphs $\G = \cXa(m,n,k,\ell)$ from item (ii) of the theorem. In fact, since the mapping $\rho$ defined as in~\eqref{eq:rho_alt} clearly is an automorphism of $\G$ preserving $\cC$, Proposition~\ref{pro:alt_gamma2} implies that we only need to show that there exists an automorphism of $\G$ preserving $\cC$ and mapping $\vv{0,0}$ to $\vv{0,1}$. Let $\beta$ be the permutation of the vertex set of $\G$ defined by the rule
\begin{equation}
\label{eq:alt_beta}
	\beta(\vv{i,j}) = \left\{\begin{array}{ccl}
		\vv{0,1-j} & : & i = 0\\
		\vv{m-i, 1-j-\ell} & : & i \neq 0,\end{array}\right.\quad i \in \ZZ_m,\ j \in \ZZ_n.
\end{equation}
Observe that $\beta$ is well defined and is a permutation of the vertex set of $\G$ (in fact, it is an involution). Since $m$ is even, $i$ and $m-i$ are of the same parity, and so the non-links of $\G$ are clearly mapped to edges of $\G$. To see that the links are also mapped to edges of $\G$ it suffices to note that $1-j-\ell$ and $j$, as well as $i$ and $m-i-1$ are of different parity. We leave the details to the reader. 
\end{proof}

\section{Additional automorphisms}
\label{sec:alt_extraAut}

In the last part of the paper we determine those graphs from Theorem~\ref{the:altVT} that admit additional automorphisms which do not preserve the $2$-factor $\cC$. We start this rather long and tedious analysis by recording some additional isomorphisms between the $\cXa(m,n,k,\ell)$ graphs and by showing that we can restrict ourselves to the examples with $n \geq 14$.

\begin{lemma}
\label{le:alt_iso2}
Let $\G = \cXa(m,n,k,\ell)$ where the parameters $m$, $n$, $k$ and $\ell$ satisfy the conditions of Theorem~\ref{the:altVT}~(ii). Then $\G = \cXa(m,n,-k,\ell)$ and $\G \cong \cXa(m,n,k,-\ell)$. Moreover, if $n$ is divisible by $4$, then $\G \cong \cXa(m,n,k+n',\ell)$ or $\G \cong \cXa(m,n,k+n',\ell+n')$, where $n = 2n'$, depending on whether $m$ is divisible by $4$ or not, respectively.
\end{lemma}

\begin{proof}
The first claim follows from the definition of the $\cXa(m,n,k,\ell)$ graphs, while an isomorphism showing that $\G \cong \cXa(m,n,k,-\ell)$ can be obtained by mapping each vertex $\vv{i,j}$ of $\G$ to the vertex $\vv{i,-j}$ of $\cXa(m,n,k,-\ell)$. To prove the last claim suppose $n$ is divisible by $4$. Then $n'$ is even and $k' = k+n'$ is coprime to $n$. Moreover, $2k' = 2k$ and $k'^2 = k^2$, and so $k'$ satisfies all the conditions of Assumption~\ref{assump_alt2}. If $m$ is divisible by $4$, we can apply Lemma~\ref{le:alt_iso}~(i) with $q = n'$ for each $i \in \{1,2,5,6,9,10,\ldots, m-3, m-2\}$ to see that $\G \cong \cXa(m,n,k',\ell)$. If however, $m \equiv 2 \pmod{4}$ then we can apply Lemma~\ref{le:alt_iso}~(i) with $q = n'$ for each $i \in \{1,2,5,6,9,10,\ldots, m-5, m-4\}$ and Lemma~\ref{le:alt_iso}~(ii) to see that $\G \cong \cXa(m,n,k',\ell+n')$ as claimed.
\end{proof}

\begin{proposition}
\label{pro:alt_n>=14}
Let $\G = \cXa(m,n,k,\ell)$ where the parameters $m$, $n$, $k$ and $\ell$ satisfy the conditions of Theorem~\ref{the:altVT}~(ii). Then $\G$ is bipartite, is of girth at least $8$ and either $n = 10$ or $n \geq 14$. Moreover, if $n = 10$ then $\G \cong \cXa(m,10,3,0)$ which is $2$-arc-regular for $m = 4$, while for all $m > 4$ the automorphism group $\Aut(\G)$ preserves the $2$-factor $\cC$ and is of order $20m$. 
\end{proposition}

\begin{proof}
That $\G$ is bipartite follows from Construction~\ref{cons:alt} since $m$ and $\ell$ are both even. Since $\gcd(k,n) = 1$ and $2k \neq \pm 2$, it follows that $n = 10$ or $n \geq 14$. Moreover, considering the vertices at distance $3$ from $\vv{1,0}$ one finds that $\G$ has no $6$-cycles (recall that $m \geq 4$ and $2k \neq \pm 2$), thus showing that $\G$ is of girth at least $8$. 

To complete the proof assume $n = 10$. Then $2k \neq \pm 2$ implies $k = \pm 3$, and so we can assume that $k = 3$ by Lemma~\ref{le:alt_iso2}. Then $k^2 = -1$, and so Theorem~\ref{the:altVT} implies that $\ell (3 \pm 1) = 0$, forcing $\ell = 0$ (recall that $\ell$ is even). That $\cXa(4,10,3,0)$ is $2$-arc-regular can be verified by {\sc Magma}~\cite{magma}. Suppose then that $m \geq 6$. One can verify that in this case there are precisely four $8$-cycles through the edge $\vv{0,0}\vv{1,0}$ - two through each of the $3$-paths $(\vv{0,0},\vv{1,0},\vv{1,3},\vv{1,6})$ and $(\vv{0,0},\vv{1,0},\vv{1,7},\vv{1,4})$. On the other hand, there are six $8$-cycles through the edge $\vv{0,0}\vv{0,1}$ - two through each of the $3$-paths $(\vv{0,0},\vv{0,1},\vv{0,2},\vv{0,3})$ and $(\vv{0,0},\vv{0,1},\vv{0,2},\vv{1,2})$ and one through each of the $3$-paths $(\vv{0,0},\vv{0,1},\vv{m-1,1},\vv{m-1,8})$ and $(\vv{0,0},\vv{0,1},\vv{m-1,1},\vv{m-1,4})$. Together with Lemma~\ref{le:alt_regular} this shows that no automorphism of $\G$ can map the link $\vv{0,0}\vv{1,0}$ to any of the non-links $\vv{0,0}\vv{0,1}$ and $\vv{0,0}\vv{0,-1}$. Consequently, $\Aut(\G)$ preserves the $2$-factor $\cC$ and $|\Aut(\G)| = 2|V(\G)| = 20m$.
\end{proof}

\begin{lemma}
\label{le:alt_aux14}
If $\G = \cXa(m,n,k,\ell)$, where the parameters $m$, $n$, $k$ and $\ell$ satisfy the conditions of Theorem~\ref{the:altVT}~(ii) and $n \geq 14$, then $2k \neq \pm 4$. 
\end{lemma}

\begin{proof}
Suppose to the contrary that $2k = \pm 4$ and note that (replacing $k$ by $-k$ if necessary) we can then in fact assume that $2k = 4$. Thus $2(k-2) = 0$, and so since $k$ is odd we must have that $k = n' + 2$ (where $n = 2n'$) and that $n$ is not divisible by $4$. But then $2k^2 = 2(n'+2)^2 = 8 \neq \pm 2$ (recall that $n \geq 14$) contradicting Theorem~\ref{the:altVT}.
\end{proof}

We now focus on $10$-cycles of the graphs $\cXa(m,n,k,\ell)$. Let $\G = \cXa(m,n,k,\ell)$ be as in Theorem~\ref{the:altVT}~(ii) with $n \geq 14$. We make an agreement that a given $10$-cycle will be identified with any of the $20$ possible presentations of the form $(v_0,v_1,\ldots , v_9)$, obtained from any one of them by cyclic rotations and reflections. Consider the presentation $(v_0,v_1,\ldots , v_9)$ of a $10$-cycle of $\G$. We assign the sequence $c_0c_1\cdots c_9$ to this presentation where for each $j \in \ZZ_{10}$ we let $c_j$ be $0$ whenever $v_jv_{j+1}$ is a link and we let $c_j$ be $1$ otherwise. The set of all sequences obtained from $c_0c_1\cdots c_9$ by cyclic rotations and reflections (which coincides with the set of all sequences obtained in the above described way from the different presentations of this cycle) will be called the {\em code} of this $10$-cycle and represented by any of the corresponding sequences. We usually abbreviate this sequence by writing it in the ``power format'' (gathering consecutive symbols of the same kind). 

Observe that 
$$
	(\vv{0,0},\vv{0,1},\vv{0,2},\vv{1,2},\vv{1,2+k},\vv{2,2+k},\vv{2,1+k},\vv{2,k},\vv{1,k},\vv{1,0})
$$
is a $10$-cycle of $\G$ with code $1^2 0101^2 010$. We say that the $10$-cycles of $\G$ with this code are of {\em type}~$0$. To determine other possible $10$-cycles observe that no $10$-cycle can have two consecutive links and that whenever there are an odd number of consecutive non-links between two links these two links have their endvertices in the same pair of sets $V_i$ and $V_{i+1}$. Since $m$ is even and $n \geq 14$, it is now easy to see that besides those of type~$0$ the only other potential $10$-cycles of $\G$ have codes  
$$
	1^4 0 1^4 0, 1^6 0 1^2 0 \ \text{and}\ 1^3 0101010,
$$
where the latter code is only possible if $m = 4$. We say that a $10$-cycle of $\G$ is of {\em type} $1$, $2$ or $3$, depending on whether its code is $1^4 0 1^4 0$, $1^6 0 1^2 0$, or $1^3 0101010$, respectively. The next result gives necessary and sufficient conditions for the existence of $10$-cycles of each of the types $0$, $1$ and $2$, and gives the number of such $10$-cycles through a given edge (those of type~3 are somewhat special, so we treat them separately).

\begin{proposition}
\label{pro:alt_10cycles}
Let $\G = \cXa(m,n,k,\ell)$ where the parameters $m$, $n$, $k$ and $\ell$ satisfy the conditions of Theorem~\ref{the:altVT}~(ii) and $n \geq 14$. Then for each of the types $0$, $1$ and $2$ of $10$-cycles of $\G$ a necessary and sufficient condition for the existence of such $10$-cycles, together with the number of $10$-cycles of this type through each of the edges $\vv{0,0}\vv{0,1}$, $\vv{0,0}\vv{0,-1}$ and $\vv{0,0}\vv{1,0}$ is as given in the following table.
{\rm 
$$
\begin{array}{c||c|c|c|c}
	\text{type} & \text{condition} & \vv{0,0}\vv{0,1} & \vv{0,0}\vv{0,-1} & \vv{0,0}\vv{1,0} \\ 
	\hline
	0 & \text{none}  & 6 & 6 & 8 \\
	1 & 4k = \pm 4 & 4 & 4 & 2  \\
	2 & 2k = \pm 6 & 8 & 8 & 4  \\
\end{array}
$$
}
\end{proposition}

\begin{proof}
We first consider the $10$-cycles of type~$0$. Due to the code of these $10$-cycles we see that for each such $10$-cycle of $\G$ there exists a unique $i \in \ZZ_m$ such that this cycle contains two consecutive edges of the cycle $C_i$, two consecutive edges of $C_{i+2}$ and two nonconsecutive edges of the cycle $C_{i+1}$. Moreover, suppose $j \in \ZZ_n$ is such that one of the two edges from $C_{i+1}$ is $\vv{i+1,j}\vv{i+1,j+1}$ or $\vv{i+1,j}\vv{i+1,j+k}$, depending on whether $i+1$ is even or odd, respectively. Since Lemma~\ref{le:alt_aux14} implies that $2k \neq \pm 4$, and consequently also $4k \neq \pm 2$ (recall that $2k^2 = \pm 2$), it then follows that for a suitable $\delta \in \{-1,1\}$ the other edge from $C_{i+1}$ on this $10$-cycle is $\vv{i+1,j+2\delta k}\vv{i+1,j+1+2\delta k}$ or $\vv{i+1,j+2\delta}\vv{i+1,j+2\delta + k}$, depending on whether $i+1$ is even or odd, respectively. It is now clear that the existence of $10$-cycles of type~$0$ is not subject to any extra condition on the parameters. It is also easy to verify that for each of the three given edges the number of $10$-cycles of type~$0$ through this edge is as stated in the above table.

That the two conditions for the $10$-cycles of types~$1$ and $2$ are as stated in the table is clear (recall that $2k^2 = \pm 2$, and so the condition $6k = \pm 2$ is equivalent to the condition $2k = \pm 6$). Since $n \geq 14$ and $\gcd(k,n) = 1$, we find that $4k \neq -4k$ and $6k \neq -6k$. It is now easy to confirm the stated numbers of $10$-cycles through each of the three given edges from the above table.
\end{proof}



We now consider the $10$-cycles of type~$3$. Recall that for $10$-cycles of type~$3$ to exist we require $m = 4$ to hold. Moreover, for each such $10$-cycle there exists an $i \in \ZZ_4$ such that this $10$-cycle has three consecutive edges of the cycle $C_i$, one edge from each of the other three cycles from $\cC$ and four links. In view of the existence of the automorphism $\gamma$ from the proof of Proposition~\ref{pro:alt_gamma2} we can assume that $i = 0$. There thus exist some $\varepsilon_0, \varepsilon_1, \varepsilon_2, \varepsilon_3 \in \{-1,1\}$ such that $3\varepsilon_0 + \varepsilon_1 k + \varepsilon_2 + \varepsilon_3 k + \ell = 0$. Theorem~\ref{the:altVT} implies that $2k \neq \pm 2$ and $\ell k = \pm \ell$, and so there are three essentially different possibilities depending on whether $\ell$ is of the form $\pm 4$, $\pm 2 \pm 2k$ or $\pm 4 \pm 2k$. We say that the $10$-cycles corresponding to these three types of conditions are of types $3.1$, $3.2$ and $3.3$, respectively. Observe that in the case that $10$-cycles of type~3.1 exist, Theorem~\ref{the:altVT} implies that $4k = \pm 4$, and so we also have that $\ell = \pm 4k$.

\begin{lemma}
\label{le:alt_t3_1}
Let $\G = \cXa(4,n,k,\ell)$ where the parameters $n$, $k$ and $\ell$ satisfy the conditions of Theorem~\ref{the:altVT}~(ii) with $m = 4$ and $n \geq 14$. Suppose $\G$ possesses $10$-cycles of type~$3.3$. Then either $\G = \cXa(4,20,k,10)$ for some $k \in \{\pm 3, \pm 7\}$ in which case $\G$ is $2$-arc-regular, or $\Aut(\G)$ is regular and thus preserves the $2$-factor $\cC$.
\end{lemma}

\begin{proof}
By Lemma~\ref{le:alt_iso2} we can assume that $\ell = 4 + 2k$. Let $\delta \in \{-1,1\}$ be such that $2k^2 = 2\delta$. Since $m = 4$, Theorem~\ref{the:altVT} implies that $\ell k = \varepsilon \ell$ for some $\varepsilon \in \{-1,1\}$. Suppose first that $\varepsilon = 1$. Then $0 = \ell(k-1) = 2k-4+2\delta$, and so Theorem~\ref{the:altVT} implies that $\delta = -1$ and $2k = 6$. By Lemma~\ref{le:alt_swap} we also have that $2\ell = 0$, and so $0 = 4k+8 = 20$ holds in $\ZZ_n$. Therefore, $n \geq 14$ implies that $n = 20$, and so $k \in \{3,13\}$ and $\ell = 10$. Using {\sc Magma} one can verify that the graph $\cXa(4,20,3,10)$ (which is isomorphic to $\cXa(4,20,13,10)$ by Lemma~\ref{le:alt_iso2}) is $2$-arc-regular.

Suppose now that $\varepsilon = -1$. Then
$$
	0 = \ell(k+1) = 2k^2 + 6k + 4 = 6k + 4 + 2\delta.
$$
If $\delta = -1$ (in which case $2\ell = 0$), we thus have that $6k + 2 = 0 = 4k + 8$, forcing $2k = 6$. Then $4k = 12$, and so again $n = 20$. As above, $k \in \{3,13\}$ and $\ell = 10$ holds. We are thus left with the possibility that $\delta = 1$, in which case $6k = -6$. If $\G$ has $10$-cycles of type~$1$, then $4k = 4$ (recall that $2k \neq -2$). In this case $0 = 12k+12 = 24$, and so in fact $n = 24$ (recall that $n \geq 14$). Then $4k = 4$ and $6k = -6$ imply that $k \in \{7,19\}$ and consequently $\ell = 18$. Using {\sc Magma} one can verify that the graph $\cXa(4,24,7,18)$ (which is isomorphic to $\cXa(4,24,19,18)$ by Lemma~\ref{le:alt_iso2}) has a regular automorphism group. 

We can thus assume that $\G$ has no $10$-cycles of type~$1$. Since $6k = -6$ and $\gcd(k,n) = 1$, Proposition~\ref{pro:alt_10cycles} implies that there are also no $10$-cycles of type~$2$. Since there are no $10$-cycles of type~$1$, Theorem~\ref{the:altVT} implies that there are also no $10$-cycles of type~$3.1$. Moreover, if $\ell = \pm 2 \pm 2k$, then Lemma~\ref{le:alt_aux14} implies that $-2-2k = \ell = 4 + 2k$, contradicting $-6 = 6k$. This shows that we only have the $10$-cycles of types~$0$ and~$3.3$. Moreover, $\ell \notin \{4 - 2k, -4 + 2k, -4-2k\}$ (recall that $2k \neq 2$), and so $\ell = 4 + 2k$ is the only condition of type $\ell = \pm 4 \pm 2k$ giving rise to $10$-cycles of type~$3.1$ (note however that $6k = -6$ implies that $4k+2 = 6k-2k+2 = -2k-4 = -\ell$, and so $\ell = -2 - 4k$ also holds). It is now not difficult to see that there are precisely three $10$-cycles of type~$3$ through the edge $\vv{0,0}\vv{0,1}$, namely
\begin{equation}
\label{eq:alt_10cyc_1}
\begin{array}{c}
	(\vv{0,1},\vv{0,0},\vv{0,-1},\vv{0,-2},\vv{1,-2},\vv{1,-2-k},\vv{2,-2-k},\vv{2,-3-k},\vv{3,-3-k},\vv{3,-3-2k}), \\
	(\vv{0,1},\vv{0,0},\vv{1,0},\vv{1,-k},\vv{2,-k},\vv{2,-1-k},\vv{2,-2-k},\vv{2,-3-k},\vv{3,-3-k},\vv{3,-3-2k}), \\
	(\vv{0,1},\vv{0,0},\vv{1,0},\vv{1,-k},\vv{2,-k},\vv{2,-1-k},\vv{3,-1-k},\vv{3,-1-2k},\vv{0,3},\vv{0,2})\\
\end{array}
\end{equation}
and precisely three through the edge $\vv{0,0}\vv{0,-1}$, namely the first from~\eqref{eq:alt_10cyc_1} and 
\begin{equation}
\label{eq:alt_10cyc_2}
\begin{array}{c}
	(\vv{0,-1},\vv{0,0},\vv{1,0},\vv{1,k},\vv{1,2k},\vv{1,3k},\vv{2,3k},\vv{2,1+3k},\vv{3,1+3k},\vv{3,1+4k}), \\
	(\vv{0,-1},\vv{0,0},\vv{1,0},\vv{1,k},\vv{2,k},\vv{2,1+k},\vv{3,1+k},\vv{3,1+2k},\vv{3,1+3k},\vv{3,1+4k}). \\
\end{array}
\end{equation}
Moreover, the last two $10$-cycles from~\eqref{eq:alt_10cyc_1} and the two $10$-cycles from \eqref{eq:alt_10cyc_2} are the only $10$-cycles of type~$3$ through the edge $\vv{0,0}\vv{1,0}$. Therefore, Proposition~\ref{pro:alt_10cycles} implies that the $2$-factor $\cC$ is $\Aut(\G)$-invariant, and consequently Lemma~\ref{le:alt_regular} shows that $\Aut(\G)$ is regular (recall that $2\ell \neq 0$).
\end{proof}

\begin{proposition}
\label{pro:alt_4kneq4}
Let $\G = \cXa(m,n,k,\ell)$ where the parameters $m$, $n$, $k$ and $\ell$ satisfy the conditions of Theorem~\ref{the:altVT}~(ii) and $n \geq 14$. If the $2$-factor $\cC$ is not preserved by $\Aut(\G)$, then either $4k = \pm 4$ or $\G = \cXa(4,20,k,10)$ for some $k \in \{\pm 3, \pm 7\}$.
\end{proposition}

\begin{proof}
Suppose that the $2$-factor $\cC$ is not preserved by $\Aut(\G)$ and that $4k \neq \pm 4$. Proposition~\ref{pro:alt_10cycles} then implies that $\G$ must have $10$-cycles of type~$3$, and so $m = 4$. By Lemma~\ref{le:alt_t3_1} it suffices to show that we must have $10$-cycles of type~$3.3$. By way of contradiction suppose that this is not the case. Since $4k \neq 4$, Theorem~\ref{the:altVT} implies that there are no $10$-cycles of type~$3.1$, and so the only $10$-cycles of type~$3$ are those of type~$3.2$. By Lemma~\ref{le:alt_iso2} we can thus assume that $\ell = 2 + 2k$. Since $4k \neq \pm 4$, we have that $2\ell \neq 0$, and so Lemma~\ref{le:alt_swap} implies that $2k^2 = 2$. By Theorem~\ref{the:altVT} there is some $\varepsilon \in \{-1,1\}$ such that 
$$
	0 = \ell(k+\varepsilon) = 2k^2 + 2(1+\varepsilon)k + 2\varepsilon = 2(1+\varepsilon)k + 2(1+\varepsilon),
$$
forcing $\varepsilon = -1$. If $\G$ has $10$-cycles of type~$2$, then Proposition~\ref{pro:alt_10cycles} implies that $36 = (2k)^2 = 4k^2 = 4$, and so $n \in \{16,32\}$. However, since $2k = \pm 6$ forces $2 = 2k^2 = \pm 18$, we find that $n = 16$ and $4k = \pm 4$, a contradiction.

We are thus left with the possibility that the only $10$-cycles of $\G$ are those of types~$0$ and~$3.2$ (which correspond only to the condition $\ell = 2+2k$). It is now easy to verify that each of the edges $\vv{0,0}\vv{0,1}$ and $\vv{0,0}\vv{1,0}$ lies on four different $10$-cycles of type~$3$, while the edge $\vv{0,0}\vv{0,-1}$ lies on just two. But then Proposition~\ref{pro:alt_10cycles} implies that the $2$-factor $\cC$ is $\Aut(\G)$-invariant, a contradiction.
\end{proof}

We finally analyze the graphs $\cXa(m,n,k,\ell)$ from Theorem~\ref{the:altVT}~(ii) that have $10$-cycles of type~$1$. There are two essentially different possibilities depending on whether the graph has only $10$-cycles of types~$0$ and $1$ or not. Note that by Proposition~\ref{pro:alt_10cycles} in the latter case $10$-cycles of type~$3$ must exist if the $2$-factor $\cC$ is not invariant under the full automorphism group of the graph. We first analyze this situation.

\begin{proposition} 
\label{pro:alt_type1_and_3}
Let $\G = \cXa(m,n,k,\ell)$ where the parameters $m$, $n$, $k$ and $\ell$ satisfy the conditions of Theorem~\ref{the:altVT}~(ii) and where $n \geq 14$ and $4k = \pm 4$. If the $2$-factor $\cC$ is not preserved by $\Aut(\G)$ and $\G$ possesses $10$-cycles of type~$3$, then there exists a positive integer $n_0$ such that up to isomorphisms from Lemma~\ref{le:alt_iso2} the graph $\G$ is one of $\cXa(4,16n_0,4n_0+1,8n_0+4)$ and $\cXa(4,32n_0,8n_0+1,4)$.
\end{proposition}

\begin{proof}
By Lemma~\ref{le:alt_iso2} we can assume that $4k = 4$. Since $2k \neq 2$, $n$ is divisible by $4$, say $n = 4n_1$ for some $n_1 \geq 4$, and $2(k-1) = 2n_1$. By Lemma~\ref{le:alt_iso2} we can replace $k$ by $k + 2n_1$ if necessary (since $m = 4$), and so we can in fact assume that $k = n_1+1$. Therefore, $8 \mid n$ (recall that $k$ is odd) and $2k^2 = 2$. 

Suppose that $\cC$ is not $\Aut(\G)$-invariant and that $\G$ possesses $10$-cycles of type~$3$. By Lemma~\ref{le:alt_t3_1} these $10$-cycles cannot be of type~$3.3$. Since Lemma~\ref{le:alt_iso2} allows us to replace $\ell$ by $-\ell$ if necessary we can thus assume that $\ell \in \{4, 2k-2, 2k+2\} = \{4, 2n_1, 2n_1+4\}$. If $n = 16$, then $\ell \in \{4,8,12\}$. Using {\sc Magma} we can verify that the graph $\cXa(4,16,5,8)$ has $128$ automorphisms, so that by Lemma~\ref{le:alt_regular} the $2$-factor $\cC$ is $\Aut(\G)$-invariant, while the graph $\cXa(4,16,5,4)$ (which is isomorphic to $\cXa(4,16,5,12)$), has $256$ automorphisms, and so $\cC$ is not $\Aut(\G)$-invariant. 

For the rest of the proof we can thus assume that $n > 16$. Since in this case $2k = 2n_1+2 \neq \pm 6$, Proposition~\ref{pro:alt_10cycles} implies that we do not have $10$-cycles of type~$2$ in $\G$. Then $4$, $2k-2$, $2k+2$ and $-2k-2$ are four different elements of $\ZZ_n$ (note however that $4 = 4k$ and $2k-2 = -2k+2$), and so it is now not difficult to determine the $10$-cycles of type~$3$ through each of the three edges incident with $\vv{0,0}$, depending on whether $\ell$ is $4$, $2k-2$ or $2k+2$. We leave it to the reader to verify that the number of $10$-cycles of type~$3$ through the corresponding edge (depending on the choice for $\ell$) is as given in the following table.
$$
\begin{array}{c||c|c|c}
	\ell & \vv{0,0}\vv{0,1} & \vv{0,0}\vv{0,-1} & \vv{0,0}\vv{1,0} \\ 
	\hline
	4 & 8 & 4 & 8 \\
	2k-2 & 6 & 6 & 8  \\
	2k+2 & 4 & 2 & 4  \\
\end{array}
$$
Since $\cC$ is not $\Aut(\G)$-invariant, Proposition~\ref{pro:alt_10cycles} implies that $\ell \in \{4, 2k+2\}$ and that the edges $\vv{0,0}\vv{0,1}$ and $\vv{0,0}\vv{1,0}$ are in the same $\Aut(\G)$-orbit, while $\vv{0,0}\vv{0,-1}$ is in a different $\Aut(\G)$-orbit. We denote the latter orbit by $\O$. Let $\rho \in \Aut(\G)$ be as in~\eqref{eq:rho_alt} and $\gamma \in \Aut(\G)$ be as in the proof of Proposition~\ref{pro:alt_gamma2} and note that the $\varepsilon$ from that proof must be equal to $1$ (since that proof shows that $\ell k = \varepsilon \ell$). The action of the subgroup $\la \rho, \gamma\ra$ then reveals that
\begin{equation}
\label{eq:alt_orbit}
	\O = \{\vv{i,2j-1}\vv{i,2j} \colon i \in \{0,2\},\ j \in \ZZ_n\} \cup \{\vv{i,2j}\vv{i,2j+k} \colon i \in \{1,3\},\ j \in \ZZ_n\}.
\end{equation}
By Theorem~\ref{the:altVT} there exists a vertex-transitive subgroup of $\Aut(\G)$ preserving $\cC$, and so as $\O$ is an $\Aut(\G)$-orbit and $\cC$ is not $\Aut(\G)$-invariant, there exists some $\theta \in \Aut(\G)$ fixing $\vv{0,0}$ but interchanging $\vv{0,1}$ with $\vv{1,0}$. Of course, $\vv{0,-1}$ is fixed by $\theta$. We claim that we can assume that $\theta(\vv{0,-2}) = \vv{3,-1-\ell}$. If this is not the case, then $\theta$ fixes each of $\vv{0,-2}$ and $\vv{0,-3}$ (recall that $\vv{0,-2}\vv{0,-3} \in \O$). Note that there is a unique $10$-cycle of $\G$ through the path $(\vv{0,-4},\vv{0,-3},\vv{0,-2},\vv{0,-1},\vv{0,0})$ and that it is of type~$1$ and continues through the vertex $\vv{1,0}$, which is not fixed by $\theta$. It thus follows that $\theta(\vv{0,-4}) = \vv{3,-3-\ell}$. It is now easy to see that the automorphism $\theta\rho\theta\rho^{-1}$ fixes each of $\vv{0,0}$ and $\vv{0,-1}$ but none of $\vv{0,1}$ and $\vv{0,-2}$. There thus exists a $\theta \in \Aut(\G)$ such that 
$$
	\theta(\vv{0,-2}) = \vv{3,-1-\ell},\ \theta(\vv{0,-1}) = \vv{0,-1},\ \theta(\vv{0,0}) = \vv{0,0}\ \text{and}\ \theta(\vv{0,1}) = \vv{1,0}.
$$
By~\eqref{eq:alt_orbit} we then have that $\theta(\vv{0,-3}) = \vv{3,-1-k-\ell}$ and $\theta(\vv{0,2}) = \vv{1,k}$. Note that there is a unique $10$-cycle through the path $(\vv{0,-2},\vv{0,-1},\vv{0,0},\vv{0,1},\vv{0,2})$ and that the next vertex on it is $\vv{1,2}$. Since this path is mapped by $\theta$ to the path $(\vv{3,-1-\ell},\vv{0,-1},\vv{0,0},\vv{1,0},\vv{1,k})$, which clearly lies on a $10$-cycle of type~$0$, it thus follows that $\theta(\vv{1,2}) = \vv{1,2k}$. Therefore, $\theta(\vv{0,3}) = \vv{2,k}$ and $\theta(\vv{0,4}) = \vv{2,k+1}$. An inductive approach now shows that each path of length $4$ of the form $(\vv{0,2j},\vv{0,2j+1},\vv{0,2j+2},\vv{0,2j+3},\vv{0,2j+4})$, which lies on a unique $10$-cycle of $\G$ (which is of type~$1$) is mapped by $\theta$ to a path of length $4$ that lies on a $10$-cycle of type~$0$, and that all edges of the form $\vv{0,2j}\vv{0,2j+1}$ are mapped to links of $\G$. It is now easy to determine the action of $\theta$ on $V_0$. In particular, $\theta(\vv{0,j}) \in V_0$ if and only if one of $j$ and $j+1$ is divisible by $8$, while $\theta(\vv{0,j}) \in V_1$ if and only if one of $j-1$ and $j-2$ is divisible by $8$. Moreover, since $\theta^2$ fixes the path $(\vv{0,-2},\vv{0,-1},\vv{0,0},\vv{0,1},\vv{0,2})$, it must also fix the unique $10$-cycle through it pointwise, and so also fixes $\vv{0,3}$ and hence $\vv{0,4}$. An inductive approach thus shows that $\theta^2$ fixes each vertex of $V_0$ pointwise and it then easily follows that $\theta$ is an involution.

To complete the proof note that $\theta(\vv{0,2}) = \vv{1,k}$ and $\theta(\vv{1,2}) = \vv{1,2k}$. Since $\vv{1,2}\vv{1,2-k} \notin \O$, $\theta$ interchanges $\vv{1,2-k}$ and $\vv{0,2k}$. But then $2k - 2 = 2n_1$ must be divisible by $8$, showing that $n$ is divisible by $16$. Finally, if $\ell = 4$, then since $2k+2+\ell = 2n_1 + 8 \neq 8$, the above inductive approach shows that $\theta$ does not fix $\vv{0,8}$ (it maps it to $\vv{0,2n_1+8}$), while it fixes each vertex of the form $\vv{0,16j}$, $j \in \ZZ_n$. Thus $2n_1 + 8$ is not divisible by $16$, showing that $n$ in fact must be divisible by $32$.
\end{proof}

In Proposition~\ref{pro:alt_swap_at_<0,0>_auto} we will show that the graphs $\cXa(4,16n_0,4n_0+1,8n_0+4)$ and $\cXa(4,32n_0,8n_0+1,4)$, $n_0 \geq 1$, indeed do admit an automorphism $\theta$ not preserving the $2$-factor $\cC$. Note that the proof of Proposition~\ref{pro:alt_type1_and_3} shows that the stabilizer of the vertex $\vv{0,0}$ in these graphs is of order at most $4$. Not to make this paper longer than it already is we mention without proof that it can in fact be shown that the stabilizer is in fact of order $2$, except for the graphs $\cXa(4,16n_0,4n_0+1,8n_0+4)$ with $n_0$ odd for which it is indeed of order $4$.  

We finally focus on the examples that only have $10$-cycles of types~$0$ and $1$. 

\begin{lemma}
\label{le:alt_type1_fix}
Let $\G = \cXa(m,n,k,\ell)$ where the parameters $m$, $n$, $k$ and $\ell$ satisfy the conditions of Theorem~\ref{the:altVT}~(ii) and where in addition $n \geq 14$, $4k = \pm 4$ and $\G$ has no $10$-cycles of types~$2$ and $3$. Then the only automorphism of $\G$ fixing a vertex and all of its neighbors is the identity. Consequently, $|\Aut(\G)| \leq 6mn$.
\end{lemma}

\begin{proof}
Suppose that $\G$ only has $10$-cycles of types~$0$ and~$1$ and let $\tau \in \Aut(\G)$ fix $\vv{0,0}$ and each of its neighbors. It is easy to see that there are precisely five $10$-cycles through the $2$-path $(\vv{0,-1},\vv{0,0},\vv{0,1})$, two of type~$0$ and three of type~$1$. Since this is an odd number, the number of these $10$-cycles that pass through $\vv{0,2}$ (which happens to be $2$) is different from the number of these $10$-cycles that pass through $\vv{m-1,1-\ell}$, and so $\tau$ fixes $\vv{0,2}$. Similarly, $\tau$ fixes $\vv{0,-2}$. There are three $10$-cycles through the $3$-path $(\vv{0,2},\vv{0,1},\vv{0,0},\vv{1,0})$, and so (as this is again an odd number) $\tau$ must fix each of $\vv{1,k}$ and $\vv{1,-k}$. Since by Theorem~\ref{the:altVT} the graph $\G$ is vertex-transitive this in fact shows that if some automorphism of $\G$ fixes a vertex and all of its three neighbors, then it also fixes each neighbor of each of these three neighbors. As $\G$ is connected, this shows that $\tau$ is the identity.
\end{proof}

\begin{proposition}
\label{pro:alt_swap_at_<0,0>}
Let $\G = \cXa(m,n,k,\ell)$ where the parameters $m$, $n$, $k$ and $\ell$ satisfy the conditions of Theorem~\ref{the:altVT}~(ii) and where in addition $n \geq 14$, $4k = 4$ and $\G$ has no $10$-cycles of types~$2$ and $3$. If $\G$ admits an automorphism fixing the vertex $\vv{0,0}$ and interchanging the vertex $\vv{1,0}$ with one of $\vv{0,1}$ and $\vv{0,-1}$, then up to isomorphisms from Lemma~\ref{le:alt_iso2} there exist integers $m_0, n_0$ and $\ell_0$, where $m_0 \geq 2$, $n_0 \geq 1$, $0 \leq \ell_0 < 4n_0$ and $8n_0 \mid (m_0n_0+\ell_0-1)(m_0n_0+\ell_0+3)$, such that $m = 2m_0$, $n = 4n_0 m$, $k = n_0m+1$ and $\ell = \ell_0 m$.
\end{proposition}

\begin{proof}
Since $4(k-1) = 0$, Theorem~\ref{the:altVT} implies that $n = 4n_1$ for some $n_1 \geq 4$ and $k \in \{n_1 + 1, 3n_1 + 1\}$. By Lemma~\ref{le:alt_iso2} we can thus assume that $k = n_1+1$ (but note that if $m \equiv 2 \pmod{4}$ we may have to replace $\ell$ by $\ell + 2n_1$). Since $k$ is odd, $n$ is divisible by $8$.

Suppose that $\theta \in \Aut(\G)$ fixes the vertex $\vv{0,0}$ and interchanges $\vv{1,0}$ with one of $\vv{0,1}$ and $\vv{0,-1}$. Lemma~\ref{le:alt_iso2} and its proof show that if $\theta(\vv{1,0}) = \vv{0,-1}$ we can simply replace $\ell$ by $-\ell$ and then the corresponding $\theta$ for the graph $\cXa(m,n,k,-\ell)$ will exchange $\vv{1,0}$ with $\vv{0,1}$. We can thus assume that $\theta$ interchanges the vertices $\vv{0,1}$ and $\vv{1,0}$. 

We determine the action of $\theta$ on $\G$ by inspecting its action on the $10$-cycles of $\G$. The only $10$-cycles of $\G$ are those of types~$0$ and~$1$, and so it is easy to see that for each $3$-path consisting of $3$ non-links there are precisely two $10$-cycles through it (both of type~$1$), for each $3$-path consisting of two consecutive non-links and a link there are precisely three $10$-cycles through it (two of type~$0$ and one of type~$1$), and for each $3$-path consisting of two links and a non-link there are precisely two $10$-cycles through it (both of type~$0$). The $3$-paths whose middle edge is a link come in two ``flavours''. Each such 3-path of course lies on two $10$-cycles of type~$0$. But some of them also lie on a $10$-cycle of type~$1$, while some do not. For instance, for $i$ even and $j \in \ZZ_n$ the $3$-path $(\vv{i,2j+1},\vv{i,2j},\vv{i+1,2j},\vv{i+1,2j+k})$ lies on a $10$-cycle of type~$1$, while $(\vv{i,2j+1},\vv{i,2j},\vv{i+1,2j},\vv{i+1,2j-k})$ does not. We say that a $3$-path is of {\em type}~$2$ or~$3$, depending on whether there are two or three $10$-cycles through it, respectively.

Now, observe first that by Lemma~\ref{le:alt_type1_fix} the automorphism $\theta$ is an involution and since it fixes both $\vv{0,0}$ and $\vv{0,-1}$, it must interchange $\vv{0,-2}$ with $\vv{m-1,-1-\ell}$. As the $3$-path $(\vv{0,-1},\vv{0,0},\vv{0,1},\vv{0,2})$ is of type~$2$ and $(\vv{0,-1},\vv{0,0},\vv{1,0},\vv{1,-k})$ is of type~$3$, $\theta$ interchanges $\vv{0,2}$ with $\vv{1,k}$. Next, since the $3$-path $(\vv{0,0},\vv{0,1},\vv{0,2},\vv{0,3})$ is of type~$2$, we see that $\theta$ interchanges $\vv{0,3}$ with $\vv{2,k}$. Continuing in this way we find that each edge of the form $\vv{0,2j}\vv{0,2j+1}$ is mapped to a link and each edge of the form $\vv{0,2j+1}\vv{0,2j+2}$ to a non-link where if $\theta(\vv{0,2j+1}) = \vv{i',j'}$ for appropriate $i'$ and $j'$ then $\theta(\vv{0,2j+2})$ is $\vv{i',j'+1}$ or $\vv{i',j'+k}$, depending on whether $i'$ is even or odd, respectively. Writing $m = 2m_0$ (recall that $m$ is even) we thus see that $\theta(\vv{0,2m-2}) = \theta(\vv{0,2(m-1)}) = \vv{m-1,(m_0-1)(k+1)+k}$, and so
$$
	\theta(\vv{0,2m-1}) = \vv{0,m_0(k+1)+\ell-1}\ \text{and}\ \theta(\vv{0,2m}) = \vv{0,m_0(k+1)+\ell}.
$$
This shows that $\theta(\vv{0,j}) \in V_0$ if and only if one of $j$ and $j+1$ is divisible by $2m$ and in particular $\gcd(m_0(k+1)+\ell,n) = 2m$. Thus $2m$ divides $n$ and (since $k+1$ is even) $m$ divides $\ell$, that is, $\ell = \ell_0 m$ for some $\ell_0 \geq 0$. Moreover, $\theta(\vv{0,j}) \in V_1$ if and only if one of $j-2$ and $j-1$ is divisible by $2m$. Since $\theta(\vv{0,1}) = \vv{1,0}$, $\theta(\vv{0,2}) = \vv{1,k}$ and $\theta(\vv{0,3}) = \vv{2,k}$, the fact that $\theta$ is an involution implies that $\theta(\vv{1,2k}) = \vv{1,2}$, and so $\theta(\vv{0,2k}) \in V_1$ (as $\theta(\vv{1,k}) = \vv{0,2} \in V_0$). Therefore, $2k-2$ is divisible by $2m$, and so $2m$ divides $2n_1$. This thus proves that $n = 4mn_0$ for some integer $n_0 \geq 1$. 

To complete the proof write $t = m_0(k+1)+\ell$, recall that $\theta(\vv{0,2m}) = \vv{0,t}$ and note that $t = m(m_0n_0+1+\ell_0)$. Since $\gcd(t,n) = 2m$, we have that $m_0n_0+\ell_0+1 = 2t_0$ for some $t_0 \geq 0$. Now, $\theta$ is an involution, and so $\vv{0,2m} = \theta(\vv{0,t}) = \theta(\vv{0,2mt_0})$, implying that
$$
	4m_0 = 2m = t_0 t = 4m_0 t_0^2 = m_0(m_0n_0+\ell_0+1)^2
$$
must hold in $\ZZ_n$. Consequently, $8n_0$ divides $(m_0n_0+\ell_0-1)(m_0n_0+\ell_0+3)$ as claimed.
\end{proof}

\begin{proposition}
\label{pro:alt_rotate_at_<0,0>}
Let $\G = \cXa(m,n,k,\ell)$ where the parameters $m$, $n$, $k$ and $\ell$ satisfy the conditions of Theorem~\ref{the:altVT}~(ii) and where in addition $n \geq 14$, $4k = 4$ and $\G$ has no $10$-cycles of types~$2$ and $3$. If $\G$ admits an automorphism fixing the vertex $\vv{0,0}$ and cyclically permuting its three neighbors, then up to isomorphisms from Lemma~\ref{le:alt_iso2} there exist odd integers $m_0, n_0$ and an even integer $\ell_0$, where $m_0 \geq 3$, $n_0 \geq 1$, $0 \leq \ell_0 < 4n_0$, $l_0+m_0n_0-1$ is divisible by $4$ and $n_0 \mid (\ell_0^2 + 3)$, such that $m = 2m_0$, $n = 4n_0 m$, $k = n_0m+1$ and $\ell = \ell_0 m$.
\end{proposition}

\begin{proof}
The proof is very similar to the proof of Proposition~\ref{pro:alt_swap_at_<0,0>} so we leave out most of the details and only indicate the main steps of the proof. As in the proof of Proposition~\ref{pro:alt_swap_at_<0,0>} we first note that $n$ is divisible by $8$ and that we can assume $k = n_1+1$, where $n = 4n_1$, and we also introduce the concepts of $3$-paths of types~$2$ and~$3$.

Write $m = 2m_0$ (recall that $m$ is even) and assume that $\eta \in \Aut(\G)$ fixes $\vv{0,0}$ and maps $\vv{0,1}$ to $\vv{1,0}$ and $\vv{1,0}$ to $\vv{0,-1}$. Again, we determine the action of $\eta$ on the vertices of the form $\vv{0,j}$ inductively by considering the $\eta$-images of $3$-paths of the form $(\vv{0,j-3},\vv{0,j-2}, \vv{0,j-1}, \vv{0,j})$, which are all of type~$2$. We find that 
$$
	\eta(\vv{0,1}) = \vv{1,0},\ \eta(\vv{0,2}) = \vv{1,-k},\ \eta(\vv{0,3}) = \vv{2,-k},\ \eta(\vv{0,4}) = \vv{2,-1-k},\ldots , \eta(\vv{0,2m}) = \vv{0,\ell -m_0(1+k)}.
$$
Therefore, $\gcd(\ell-m_0(1+k),n) = 2m$, implying that $2m$ divides $n$ and $m$ divides $\ell$, that is, $\ell = m\ell_0$ for some $\ell_0 \geq 0$. We also easily establish that $\eta(\vv{1,2k})) = \vv{1,-2}$, and then $\eta(\vv{0,2k}) = \vv{1,-2+k}$. The above inductive approach determining $\eta(\vv{0,j})$ for each $j \in \ZZ_n$ thus shows that $2k \equiv 2 \pmod {2m}$, that is, $2m$ divides $2k-2 = 2n_1$. It follows that $4m$ divides $n$, that is, $n = 4mn_0$ for some $n_0 \geq 1$. Note that this implies that $\eta(\vv{0,2m}) = \vv{0,m(\ell_0-m_0n_0-1)}$.

Similarly we find that since $\eta^2$ (which is $\eta^{-1}$ by Lemma~\ref{le:alt_type1_fix}) maps $\vv{0,-1}$ to $\vv{1,0}$ and $\vv{1,0}$ to $\vv{0,1}$, it maps $\vv{0,-2m}$ to $\vv{0,\ell+m_0(1+k)}$, implying that $\gcd(\ell+m_0(1+k),n) = 2m$. Therefore, 
$$
	\gcd(m(\ell_0+m_0n_0+1),4mn_0) = 2m = \gcd(m(\ell_0-m_0n_0-1),4mn_0),
$$
implying that $\ell_0$ is even, $m_0n_0$ is odd and that $4$ divides $\ell_0+m_0n_0-1$. 

Finally, since $\eta^2(\vv{0,-2m}) = \vv{0,m_0(1+k) + \ell}$, we have that $\eta(\vv{0, m_0(k+1)+\ell}) = \vv{0,-2m}$. But since $m_0(k+1)+\ell = m(\ell_0 + m_0n_0 + 1)$ and $\ell_0+m_0n_0+1$ is even, we also find that $\eta(\vv{0,m_0(k+1)+\ell}) = \vv{0,s}$, where $s = m_0(\ell_0-m_0n_0-1)(\ell_0+m_0n_0+1)$. Therefore, 
$$
	2m + m_0(\ell_0-m_0n_0-1)(\ell_0+m_0n_0+1) = m_0(4 + (\ell_0^2 - (m_0n_0+1)^2))
$$
is divisible by $n$, which clearly implies that $n_0$ divides $4+\ell_0^2-1 = \ell_0^2+3$ as claimed. 
\end{proof}

In the following two propositions we show that the graphs from Propositions~\ref{pro:alt_type1_and_3},~\ref{pro:alt_swap_at_<0,0>} and \ref{pro:alt_rotate_at_<0,0>} indeed admit automorphisms mapping a link to a non-link. The verification that the given mappings are indeed automorphisms of the corresponding graphs is quite tedious, so we leave some details to the reader (but we do point out all of the main steps). 

\begin{proposition}
\label{pro:alt_swap_at_<0,0>_auto}
Let $m_0, n_0, \ell_0$ be nonnegative integers with $m_0 \geq 2$, $n_0 \geq 1$ and $0 \leq \ell_0 < 4n_0$ such that $8n_0 \mid (m_0n_0+\ell_0-1)(m_0n_0+\ell_0+3)$. Then the graph $\cXa(2m_0,8m_0n_0,2m_0n_0+1, 2m_0\ell_0)$ admits an automorphism fixing the vertex $\vv{0,0}$ and interchanging the vertices $\vv{0,1}$ and $\vv{1,0}$.
\end{proposition}

\begin{proof}
Let $m = 2m_0$, $n = 4mn_0$, $k = mn_0+1$, $\ell = m\ell_0$ and $\G = \cXa(m,n,k,\ell)$. Furthermore, let $t = m(m_0n_0+\ell_0+1)$. Since $m_0n_0+\ell_0-1$ and $m_0n_0+\ell_0+3$ have the same remainder modulo $4$, the assumption that $8n_0$ divides their product implies that $m_0n_0+\ell_0+1 \equiv 2 \pmod{4}$ and is divisible by no odd prime divisor of $n_0$. Thus $m_0n_0+\ell_0+1 \equiv 2t_0 \pmod{4n_0}$ for some odd $t_0$ coprime to $n_0$ with $0 < t_0 < 2n_0$. Moreover, $\gcd(t,n) = 2m$. We now define a certain mapping $\theta$ on the vertex set of $\G$ fixing the vertex $\vv{0,0}$ and interchanging the vertices $\vv{0,1}$ and $\vv{1,0}$ and then show that it is in fact an automorphism of $\G$. 

Note that for each $i$ with $0 \leq i < m_0$ each element of $\ZZ_n$ can uniquely be written in the form $-1 + i(k+1) + 2am + 4j + r$, where $0 \leq a < 2n_0$, $0 \leq j < m_0$ and $0 \leq r < 4$. Moreover, since $k$ is coprime to $n$, each element of $\ZZ_n$ can also uniquely be written in the form $i(k+1) + 2am + 4j + rk$, where again $0 \leq a < 2n_0$, $0 \leq j < m_0$ and $0 \leq r < 4$. We then set
\begin{equation}
\label{eq:alt_def_theta_1}
\theta(\vv{2i,-1+i(k+1)+2am + 4j+r}) = 
	\left\{\begin{array}{lcl}
	\vv{2j,4i+at+j(k+1)+r-1} & : & r \in \{0,1\}\\
	\vv{2j+1,4i+at+j(k+1)+(r-2)k} & : & r \in \{2,3\},
	\end{array}\right.
\end{equation}
\begin{equation}
\label{eq:alt_def_theta_2}
\theta(\vv{2i+1,i(k+1)+2am + 4j+rk}) = 
	\left\{\begin{array}{lcl}
	\vv{2j,4i+at+j(k+1)+r+1} & : & r \in \{0,1\}\\
	\vv{2j+1,4i+at+j(k+1)+2+(r-2)k} & : & r \in \{2,3\}.
	\end{array}\right.
\end{equation}

Since $4i+at$ is divisible by $4$ and $k$ is coprime to $n$, it is clear that no $\theta$-image from $\eqref{eq:alt_def_theta_1}$ can be equal to a $\theta$-image from~\eqref{eq:alt_def_theta_2}. Since $\gcd(t,n) = 2m$ and $i < m_0$, it is now clear that $\theta$ is injective and is thus a permutation of the vertex set of $\G$. 

That the non-links of $\G$ are mapped to edges of $\G$ is clear from~\eqref{eq:alt_def_theta_1} and~\eqref{eq:alt_def_theta_2}, except perhaps for the non-links where one of the vertices has $r = 3$ and the other $r = 0$. We consider the possibility that the two vertices are both in some $V_{2i}$ and leave the one where they are both in some $V_{2i+1}$ to the reader (here $4k = 4$ should be used). Let $v$ be a vertex of the form $\vv{2i,-1+i(k+1)+2am+4j+3}$. By~\eqref{eq:alt_def_theta_1} we then have that
$$
	\theta(v) = \vv{2j+1,4i+at+j(k+1)+k}.
$$ 
Observe that since $2j+1$ and $4i+at+j(k+1)+k$ are both odd, the outside neighbor of $\theta(v)$ is in $V_{2j+2}$.
We consider the $\theta$-image of the neighbor $w$ of $v$ in $V_{2i}$ whose second index equals $-1+i(k+1)+2am+4j+4$. If $j < m_0-1$, we can write 
$$
	-1+i(k+1)+2am+4j+4 = -1+i(k+1)+2am+4(j+1),
$$
and so $\theta(w) = \vv{2(j+1),4i+at+(j+1)(k+1)-1}$, which is indeed a neighbor of $\theta(v)$. If however $j = m_0-1$, then we can write 
$$
	-1+i(k+1)+2am+4j+4 = -1+i(k+1)+2(a+1)m
$$ 
(with the understanding that $a+1 = 0$ if $a = 2n_0-1$), and so $\theta(w) = \vv{0,4i+(a+1)t-1}$. Since $t$ is even, the outside neighbor of $\theta(w)$ is in $V_{m-1}$ and its second index is
$$
	4i + at + t - 1 - \ell = 4i + at + m(m_0n_0+1)-1 = 4i + at + m_0(k+1)-1,
$$
which is precisely the second index of $\theta(v)$.

To complete the proof we need to verify that the links of $\G$ are also mapped to edges of $\G$. Consider first a link of the form
\begin{equation}
\label{eq:alt_proof_theta_1}
	\vv{2i,-1+i(k+1)+2am+4j+r}\vv{2i+1,-1+i(k+1)+2am+4j+r}
\end{equation}
and note that in this case $r \in \{1,3\}$. That this link is mapped to an edge of $\G$ in the case of $r = 1$ follows directly from~\eqref{eq:alt_def_theta_1} and~\eqref{eq:alt_def_theta_2}. If $r = 3$, then the first of these two vertices is mapped to $\vv{2j+1,4i+at+j(k+1)+k}$. Note that $2k = 2n_0m+2$ and write
$$
	-1+i(k+1)+2am+4j+3 = i(k+1)+2(a-n_0)m+4j+2k,
$$
where we replace $a-n_0$ by $a+n_0$ if $a < n_0$ (note that $4n_0m = 0$ in $\ZZ_n$). The second vertex in~\eqref{eq:alt_proof_theta_1} is thus mapped by $\theta$ to $\vv{2j+1,4i+(a-n_0)t+j(k+1)+2}$. Since $\gcd(t,n) = 2m$ we have that $n_0t = 4m_0n_0 = 2k-2 = -2k+2$, and so the link from~\eqref{eq:alt_proof_theta_1} is mapped to an edge of $\G$. 

That the links of the form $\vv{2i+1,i(k+1)+2am+4j+rk}\vv{2i+2,i(k+1)+2am+4j+rk}$, where $i < m_0-1$ (and $r \in \{1,3\}$) are mapped to edges of $\G$ is verified in a similar way as was done in the previous paragraph, and so we leave this to the reader. We finally consider the links of the form 
\begin{equation}
\label{eq:alt_proof_theta_2}
	\vv{m-1,(m_0-1)(k+1)+2am+4j+rk}\vv{0,(m_0-1)(k+1)+2am+4j+rk+\ell},
\end{equation}
where $r \in \{1,3\}$. Observe first that $m_0(k+1)+\ell = t = 2t_0m$, and so the second index of this vertex from $V_0$ can be written as 
$$
	2t_0m + 2am + 4j + (r-1)k - 1 = -1 + 2(t_0+a)m + 4j + (r-1)k.
$$
Moreover, by assumption
$$
	t_0t - 4m_0 = m_0((2t_0)^2 - 4) = m_0((m_0n_0+\ell_0+1)^2-4) = 0,
$$
that is, $t_0t = 4m_0$. For $r = 1$ the first of the vertices from~\eqref{eq:alt_proof_theta_2} is mapped to $\vv{2j,4(m_0-1)+at+j(k+1)+2}$, while the second is mapped to 
$$
	\vv{2j,(t_0+a)t+j(k+1)-1} = \vv{2j,4m_0+at+j(k+1)-1},
$$
which is indeed a neighbor of $\vv{2j,4(m_0-1)+at+j(k+1)+2}$. The situation with $r = 3$ is very similar and is left to the reader. 
This finally proves that $\theta$ is an automorphism of $\G$. That it fixes $\vv{0,0}$ and interchanges $\vv{0,1}$ with $\vv{1,0}$ follows from~\eqref{eq:alt_def_theta_1} and~\eqref{eq:alt_def_theta_2}.
\end{proof}

\begin{proposition}
\label{pro:alt_rotate_at_<0,0>_auto}
Let $m_0, n_0$ be odd positive integers with $m_0 \geq 3$ and $n_0 \geq 1$, and let $\ell_0$ be an even integer with $0 \leq \ell_0 < 4n_0$ such that $m_0n_0+\ell_0-1$ is divisible by $4$ and $n_0$ divides $\ell_0^2+3$. Then the graph $\cXa(2m_0,8m_0n_0,2m_0n_0+1, 2m_0\ell_0)$ admits an automorphism fixing the vertex $\vv{0,0}$ and cyclically permuting its three neighbors.
\end{proposition}

\begin{proof}
The proof is very similar to the proof of Proposition~\ref{pro:alt_swap_at_<0,0>_auto}, which is why we leave out most of the details and only point out the slight differences. Again, we let $m = 2m_0$, $n = 4mn_0$, $k = mn_0+1$, $\ell = m\ell_0$ and $\G = \cXa(m,n,k,\ell)$. This time we let $t = m(\ell_0-m_0n_0-1)$. Since $\ell_0$ is even and $m_0n_0+\ell_0-1$ is divisible by $4$, it clearly follows that $\ell_0-m_0n_0-1 \equiv 2 \pmod{4}$, and so the assumption that $n_0$ is odd and divides $\ell_0^2+3$ implies that $n_0$ is coprime to $\ell_0-m_0n_0-1$. There thus exists an odd $t_0$ coprime to $n_0$ with $0 < t_0 < 2n_0$ such that $\ell_0-m_0n_0-1 \equiv 2t_0 \pmod{4n_0}$. In particular, $t = 2mt_0$ (in $\ZZ_n$) and $\gcd(t,n) = 2m$. We now define a mapping $\eta$ on $V(\G)$ fixing the vertex $\vv{0,0}$ and cyclically permuting its three neighbors and show that $\eta \in \Aut(\G)$.

Expressing each element of $\ZZ_n$ in the same form as in the proof of Proposition~\ref{pro:alt_swap_at_<0,0>_auto} we set
\begin{equation}
\label{eq:alt_def_eta_1}
\eta(\vv{2i,-1+i(k+1)+2am + 4j+r}) = 
	\left\{\begin{array}{lcl}
	\vv{2j,at-4i-j(k+1)+1-r} & : & r \in \{0,1\}\\
	\vv{2j+1,at-4i-j(k+1)+(2-r)k} & : & r \in \{2,3\},
	\end{array}\right.
\end{equation}
\begin{equation}
\label{eq:alt_def_eta_2}
\eta(\vv{2i+1,i(k+1)+2am + 4j+rk}) = 
	\left\{\begin{array}{lcl}
	\vv{2j,at-4i-j(k+1)-1-r} & : & r \in \{0,1\}\\
	\vv{2j+1,at-4i-j(k+1)-2+(2-r)k} & : & r \in \{2,3\}.
	\end{array}\right.
\end{equation}
The proof that $\eta$ is in fact a permutation of $\G$ and that it maps all non-links of $\G$ to edges of $\G$ is very similar to the corresponding proof for $\theta$ in the proof of Proposition~\ref{pro:alt_swap_at_<0,0>_auto}. The same holds for the proof that the links of $\G$ are mapped to edges of $\G$, except possibly for the links connecting the vertices from $V_{m-1}$ and $V_0$. We thus only show how this part of the proof can be carried out. 

Observe first that the assumptions that $\ell_0$ is even and $m_0n_0+\ell_0-1$ is divisible by $4$ imply that $\ell_0$ and $m_0n_0+1$ are both even but precisely one of them is divisible by $4$. Therefore, $\ell_0^2 - (m_0n_0+1)^2+4$ is divisible by $8$. Moreover, since $n_0$ divides $\ell_0^2+3$, it divides $\ell_0^2 - (m_0n_0+1)^2+4$, which thus finally shows that 
$$
	m_0(\ell_0-m_0n_0-1)(\ell_0+m_0n_0+1) = m_0(\ell_0^2 - (m_0n_0+1)^2) = -4m_0.
$$
Recall that $\ell_0-m_0n_0-1 \equiv 2 \pmod{4}$. Since $m_0n_0$ is odd, this implies that $2m_0^2n_0(\ell_0-m_0n_0-1) = 4m_0n_0$ (in $\ZZ_n$), and so 
$$
	t(1+t_0) = 2t_0m_0(2t_0+2) = m_0(\ell_0-m_0n_0-1)(\ell_0-m_0n_0+1) = 4m_0n_0-4m_0.
$$
Consider now a link of the form 
\begin{equation}
\label{eq:alt_proof_eta_1}
	\vv{m-1,(m_0-1)(k+1)+2am+4j+rk}\vv{0,(m_0-1)(k+1)+2am+4j+rk+\ell},
\end{equation}
where $r \in \{1,3\}$. If $r = 1$ then the first of these two vertices is mapped by $\eta$ to $\vv{2j,at-4m_0-j(k+1)+2}$. For the second vertex we first note that $m_0(k+1) = m(m_0n_0+1)$ and $\ell = t+m(m_0n_0+1)$, showing that we can write its second index as
$$
	-1 + 2m(m_0n_0+1) + t + 2am + 4j = -1 + 2m(m_0n_0+t_0+a+1) + 4j.
$$
Therefore, the second vertex from~\eqref{eq:alt_proof_eta_1} is mapped by $\eta$ to the vertex in $V_{2j}$ whose second index is
$$
	(m_0n_0+t_0+a+1)t - j(k+1) + 1 = - 4m_0 + at - j(k+1) + 1,
$$
thus showing that the link from~\eqref{eq:alt_proof_eta_1} is indeed mapped by $\eta$ to an edge of $\G$. The argument for $r = 3$ is similar and is left to the reader.
\end{proof}

We can now state our main result of this section.

\begin{theorem}
\label{the:alt_main}
Let $\G$ be a connected cubic vertex-transitive graph admitting a partition of its edge set into a $2$-factor $\cC$ and a $1$-factor $\cI$ such that $\G$ is of alternating cycle quotient type with respect to $\cC$. Let $m = |\cC|$ and let $n$ be such that $\G$ is of order $mn$. Then $\cC$ is $\Aut(\G)$-invariant if and only if one of the following holds:
\begin{itemize}
\itemsep = 0pt
	\item[(i)] $\G \cong \HTG(m,n,\ell)$ for some $\ell \in \ZZ_n$ of the same parity as $m$ such that none of the following three conditions is satisfied:
	\begin{itemize}
	\itemsep = 0pt
	\item[$\bullet$] $\gcd(\ell+m, n) = 2m$ and $2mn \mid (\ell^2 + 2m\ell - 3m^2)$;
	\item[$\bullet$] $\gcd(\ell-m, n) = 2m$ and $2mn \mid (\ell^2 - 2m\ell - 3m^2)$;
	\item[$\bullet$] $\gcd(\ell+m, n) = \gcd(\ell-m, n) = 2m$ and $2mn \mid (\ell^2 + 3m^2)$.
	\end{itemize}
	\item[(ii)] $m$ is even and $\G \cong \cXa(m,n,k,\ell)$ for some even $\ell \in \ZZ_n$ and odd $k \in \ZZ_n$ with $2k \neq \pm 2$ and $2k^2 = \pm 2$, where one of the two conditions stated below is satisfied but at the same time $\G$ is not isomorphic (via the isomorphisms from Lemma~\ref{le:alt_iso2}) to $\cXa(4,10,3,0)$, $\cXa(4,20,3,10)$ or a graph from Proposition~\ref{pro:alt_swap_at_<0,0>_auto} or Proposition~\ref{pro:alt_rotate_at_<0,0>_auto}:
	\begin{itemize}
	\itemsep = 0pt
		\item[$\bullet$] $\ell k = \pm \ell$ and at least one of $k^2 = \pm 1$ and $4 \mid m$ holds, or
		\item[$\bullet$] $\ell k = n' \pm \ell$, $m \equiv 2 \pmod{4}$ and $k^2 = n' \pm 1$, where $n = 2n'$.
	\end{itemize}
\end{itemize}
\end{theorem}

\begin{proof}
By Theorem~\ref{the:altVT} the only candidates for $\G$ such that $\cC$ is $\Aut(\G)$-invariant are the graphs $\HTG(m,n,\ell)$ for some $\ell \in \ZZ_n$ of the same parity as $m$ and the graphs $\cXa(m,n,k,\ell)$ from Theorem~\ref{the:altVT}~(ii). We thus need to prove that for these graphs $\cC$ is indeed $\Aut(\G)$-invariant if and only if the conditions stated in the theorem hold. 

We first consider the graphs $\G = \HTG(m,n,\ell)$, where we can simply apply the results of~\cite{Spa22}. Since $m \geq 3$, the only examples from~\cite[Theorem~1.1]{Spa22} for which the corresponding $2$-factor $\cC$ might not be $\Aut(\G)$-invariant are the Pappus graph $\HTG(3,6,3)$ and the so-called generalized prisms $\HTG(m,4,\ell)$. However, Proposition~\ref{pro:alt_n=4} shows that for the latter ones $\cC$ is $\Aut(\G)$-invariant. Thus \cite[Theorem~1.2]{Spa22} and \cite[Lemmas~4.2--4.5]{Spa22} imply that the $\HTG(m,n,\ell)$ graphs with $m \geq 3$ for which the corresponding $2$-factor $\cC$ is not $\Aut(\G)$-invariant are precisely those satisfying any of the three conditions stated in our theorem (note that the Pappus graph $\HTG(3,6,3)$ in fact satisfies all three of them).

Let now $\G = \cXa(m,n,k,\ell)$ be as in Theorem~\ref{the:altVT}~(ii). Propositions~\ref{pro:alt_n>=14},~\ref{pro:alt_4kneq4},~\ref{pro:alt_type1_and_3},~\ref{pro:alt_swap_at_<0,0>} and~\ref{pro:alt_rotate_at_<0,0>} then imply that, up to isomorphisms from Lemma~\ref{le:alt_iso2}, the only way that $\cC$ is not $\Aut(\G)$-invariant is if $\G$ is one of $\cXa(4,10,3,0)$ and $\cXa(4,20,3,10)$ or a graph from Proposition~\ref{pro:alt_swap_at_<0,0>_auto} or Proposition~\ref{pro:alt_rotate_at_<0,0>_auto} (note that the graphs from Proposition~\ref{pro:alt_type1_and_3} satisfy the conditions of Proposition~\ref{pro:alt_swap_at_<0,0>_auto}). Propositions~\ref{pro:alt_n>=14},~\ref{pro:alt_4kneq4},~\ref{pro:alt_swap_at_<0,0>_auto} and~\ref{pro:alt_rotate_at_<0,0>_auto} imply that for these graphs the $2$-factor $\cC$ is indeed not $\Aut(\G)$-invariant.
\end{proof}

\section{Concluding remarks}
\label{sec:concluding}

In the Introduction we stated that the graphs $\cXa(m,n,k,\ell)$, which were the main objects of study in this paper, are natural generalizations of the generalized Petersen graphs and of the honeycomb toroidal graphs. Indeed, if we naturally extend the definition of these graphs in Construction~\ref{cons:alt} to allow for $m \in \{1,2\}$, the graphs $\cXa(2,n,k,0)$ correspond to the generalized Petersen graphs, while the graphs $\cXa(m,n,1,\ell)$ are the honeycomb toroidal graphs. 

Let us also indicate why our claims from the end of the Introduction hold. It is easy to see that for any $m_1 \geq 1$ setting $m = 8m_1+2$ the graphs $\cXa(m,4m,m+1,0)$  and $\cXa(m,12m,3m+1,6m)$ satisfy the assumptions from Theorem~\ref{the:altVT}, as well as those of Propositions~\ref{pro:alt_swap_at_<0,0>_auto} and~\ref{pro:alt_rotate_at_<0,0>_auto}. These results then imply that these graphs are $2$-arc-transitive and then Propositions~\ref{pro:alt_10cycles} and~\ref{pro:alt_type1_and_3} and Lemma~\ref{le:alt_type1_fix} imply that they are in fact $2$-arc-regular. A similar conclusion can be made if we set $m = 8m_1-2$ and take the graphs $\cXa(m,4m,m+1,2m)$  and $\cXa(m,12m,3m+1,0)$. It is also easy to verify that for any $m_1 \geq 1$ taking $m = 8m_1-2$ the graph $\cXa(m,28m,7m+1,12m)$ satisfies the assumptions from Theorem~\ref{the:altVT} and Proposition~\ref{pro:alt_rotate_at_<0,0>_auto}. Since $2\ell \neq 0$, Lemma~\ref{le:alt_regular} and the above results show that these graphs are arc-regular.

Finally, let us briefly comment on the graphs of bialternating cycle quotient type. These graphs will be the subject of a forthcoming paper in which the examples admitting a vertex-transitive group preserving the corresponding $2$-factor $\cC$ will be classified. However, determining for which of them the full automorphism group preserves $\cC$ at the moment seems to be a very difficult task. Namely, it turns out that unlike the graphs of alternating cycle quotient type which are all of girth at most $10$, the graphs of bialternating cycle quotient type can have girths up to $14$, and so an analogous analysis as was done in Section~\ref{sec:alt_extraAut} of this paper does not seem to be feasible. Nevertheless, it seems that like the graphs $\cXa(m,n,k,\ell)$ the examples of bialternating cycle quotient type can also have various different degrees of symmetry, and so may present a rich source of graphs for future investigations on cubic vertex-transitive graphs.

\end{document}